%% file: m7-15.tex
\documentclass{gtart_h}

\input gtmonout

\volumenumber{7}
\volumename{Proceedings of the Casson Fest}
\volumeyear{2004}
\papernumber{15}
\pagenumbers{431}{491}
\received{9 September 2003}
\revised{30 July 2004}
\accepted{1 November 2004}
\published{13 December 2004}

\usepackage{amsmath, amssymb, amscd, mathrsfs, rlepsf}

\begin{document}

\newtheorem{thm}{Theorem}[subsection]
\newtheorem{lem}[thm]{Lemma}
\newtheorem{cor}[thm]{Corollary}
\newtheorem{conj}[thm]{Conjecture}
\newtheorem{qn}[thm]{Question}
\newtheorem{prob}[thm]{Problem}
\newtheorem*{claim}{Claim}
\newtheorem*{crit}{Criterion}
\newtheorem*{mapping_class_thm}{Theorem A}
\newtheorem*{bounded_orbit_thm}{Theorem B}
\newtheorem*{Euler_class_thm}{Theorem C}
\newtheorem*{Z+Z_thm}{Theorem D}

\newtheorem{defn}[thm]{Definition}
\newtheorem{construct}[thm]{Construction}
\newtheorem{note}[thm]{Notation}

\newtheorem{rmk}[thm]{Remark}
\newtheorem{exa}[thm]{Example}

\def\C{\mathbb C}
\def\CC{\mathscr C}
\def\R{\mathbb R}
\def\Z{\mathbb Z}
\def\H{\mathbb H}
\def\P{\mathscr P}
\def\E{\mathscr E}
\def\CP{\mathbb{CP}}
\def\RP{\mathbb{RP}}
\def\homeo{\text{Homeo}}
\def\thomeo{\text{H}\til{\text{omeo}}}
\def\bhomeo{\text{BHomeo}}
\def\diffeo{\text{Diffeo}}
\def\poly{\text{Poly}}
\def\isom{\text{Isom}}
\def\stab{\text{stab}}
\def\aut{\text{Aut}}
\def\hom{\text{Hom}}
\def\MCG{\text{MCG}}
\def\SL{\text{SL}}
\def\GL{\text{GL}}
\def\PSL{\text{PSL}}
\def\aff{\text{Aff}}
\def\id{\text{Id}}
\def\inte{\text{int}}
\def\index{\text{index}}
\def\fix{\text{fix}}
\def\til{\widetilde}
\def\hat{\widehat}

\title{Circular groups, planar groups, and the Euler class}
\author{Danny Calegari}

\address{Department of Mathematics, California Institute of 
Technology\\Pasadena CA, 91125, USA}
\email{dannyc@its.caltech.edu}
%\date{11/17/2004, Version 0.15}

\begin{abstract}
We study groups of $C^1$ orientation-preserving homeo\-morph\-isms of
the plane, and pursue analogies between such groups and
circul\-arly-order\-able groups.  We show that every such group with a
bounded orbit
is circularly-orderable, and show that certain generalized braid groups
are circularly-orderable.

We also show that the Euler class of $C^\infty$ diffeomorphisms of the
plane is an {\em unbounded} class, and that any closed surface group of
genus $>1$ admits a $C^\infty$ action with arbitrary Euler class. On the
other hand, we show that $\Z \oplus \Z$ actions satisfy a homological
rigidity property: every orientation-preserving $C^1$ action of $\Z
\oplus \Z$ on the plane has trivial Euler class. This gives the complete
homological classification of surface group actions on $\R^2$ in every
degree of smoothness.
\end{abstract}

\asciiabstract{%
We study groups of C^1 orientation-preserving homeomorphisms of the
plane, and pursue analogies between such groups and
circularly-orderable groups.  We show that every such group with a
bounded orbit is circularly-orderable, and show that certain
generalized braid groups are circularly-orderable.  We also show that
the Euler class of C^infty diffeomorphisms of the plane is an
unbounded class, and that any closed surface group of genus >1 admits
a C^infty action with arbitrary Euler class. On the other hand, we
show that Z oplus Z actions satisfy a homological rigidity property:
every orientation-preserving C^1 action of Z oplus Z on the plane has
trivial Euler class. This gives the complete homological
classification of surface group actions on R^2 in every degree of
smoothness.}

\asciikeywords{Euler class, group actions, surface dynamics, braid groups,
C^1 actions}

\maketitle\vspace{-5pt}

\leftskip 25pt {\small\it This paper is dedicated to Andrew Casson, on the
occasion of his 60th\break birthday.  Happy birthday, Andrew!}\rightskip 25pt

\leftskip 0pt \section{Introduction}\rightskip 0pt

We are motivated by the following question: what kinds of countable groups
$G$ act on the plane? And for a given group $G$ known to act faithfully,
what is the best possible analytic quality for a faithful action?

This is a very general problem, and it makes sense to narrow focus in
order to draw useful conclusions. Groups can be sifted through many
different kinds of strainers: finitely presented, hyperbolic, amenable,
property T, residually finite, etc. Here we propose that ``acts on
a circle'' or ``acts on a line'' is an interesting sieve to apply to
groups $G$ which act on the plane.

The theory of group actions on {\em $1$--dimensional manifolds} is
rich and profound, and has many subtle connections with algebra, logic,
analysis, topology, ergodic theory, etc. One would hope that some of the
depth of this theory would carry across to the study of group actions
on $2$--dimensional manifolds.

The most straightforward way to establish a connection between groups
which act in $1$ and $2$ dimensions is to study when the groups acting in
either case are abstractly isomorphic. Therefore we study subgroups $G <
\homeo^+(\R^2)$, and ask under what general conditions are they isomorphic
(as abstract groups) to subgroups of $\homeo^+(S^1)$.

One reason to compare the groups $\homeo^+(S^1)$ and $\homeo^+(\R^2)$
comes from their cohomology as discrete groups. A basic theorem of Mather
and Thurston says that the cohomology of both of these groups, thought
of as discrete groups, is equal, and is equal to $\Z[e]$ where $[e]$ in
dimension $2$ is free, and is the {\em Euler class}. Thus at a classical
algebraic topological level, these groups are not easily distinguished,
and we should not be surprised if many subgroups of $\homeo^+(\R^2)$
can be naturally made to act faithfully on a circle. We establish that
countable $C^1$ groups of homeomorphisms of the plane which satisfy
a certain dynamical condition --- that they have a bounded orbit ---
are all isomorphic to subgroups of $\homeo^+(S^1)$.

On the other hand, the {\em bounded cohomology} of these groups is
dramatically different. The classical Milnor--Wood inequality says that
the Euler class on $\homeo^+(S^1)$ is a {\em bounded class}. By contrast,
the Euler class on $\homeo^+(\R^2)$ is {\em unbounded}. This was known
to be true for $C^0$ homeomorphisms; in this paper we establish that
it is also true for $C^\infty$ homeomorphisms. However, a surprising
rigidity phenomenon manifests itself: for $C^1$ actions of $\Z \oplus \Z$
we show that the Euler class must always vanish, which would be implied
by boundedness.  This is surprising for two reasons: firstly, because we
show that the Euler class can take on any value for $C^\infty$ actions
of higher genus surface groups, and secondly because the Euler class can
take on any value for $C^0$ actions of $\Z \oplus \Z$. It would be very
interesting to understand the full range of this homological rigidity.

We now turn to a more precise statement of results.

\subsection{Statement of results}

Section 2 contains background material on left-orderable and
circularly-order\-able groups. This material is all standard, and is
collected there for the convenience of the reader. The main results
there are that a countable group is {\em left-orderable}
iff it admits an injective homomorphism to $\homeo^+(\R)$, and
{\em circularly-orderable} iff it admits an injective homomorphism
to $\homeo^+(S^1)$. Thus, group actions on $1$--manifolds can be
characterized in purely algebraic terms. The expert may feel free to
skip Section 2 and move on to Section 3.

Section 3 concerns $C^1$ subgroups $G < \homeo^+(\R^2)$ with {\em
bounded orbits}. Our first main result is a generalization of a theorem
of Dehornoy \cite{Dehorn_braids} about orderability of the usual braid
groups.

For $C$ a compact, totally disconnected subset of the open unit disk $D$,
we use the notation $B_C$ to denote the group of homotopy classes of
homeomorphisms of $\overline{D}\backslash C$ to itself which are fixed
on the boundary, and $B_C'$ to denote the group of homotopy classes
of (orientation-preserving) homeomorphisms which might or might not
be fixed on the boundary. Informally, $B_C$ is the  ``braid group''
of $C$. In particular, if $C$ consists of $n$ isolated points, $B_C$
is the usual braid group on $n$ strands.

\begin{mapping_class_thm}
Let $C$ be a compact, totally disconnected subset of the open unit
disk $D$.
Then $B_C'$ is circularly-orderable, and $B_C$ is left-orderable.
\end{mapping_class_thm}

Using this theorem and the Thurston stability theorem
\cite{Thurston_stability},
we show the following:

\begin{bounded_orbit_thm}
Let $G$ be a group of orientation preserving $C^1$ homeomorphisms of
$\R^2$ with a bounded orbit. Then $G$ is circularly-orderable.
\end{bounded_orbit_thm}

Section 4 concerns the Euler class for planar actions. As intimated
above,
we show that the Euler class can take on any value for $C^\infty$ actions
of higher genus surface groups:

\begin{Euler_class_thm}
For each integer $i$, there is a $C^\infty$
action $$\rho_i\co\pi_1(S) \to \diffeo^+(\R^2)$$
where $S$ denotes the closed surface of genus $2$, satisfying
$$\rho_i^*([e])([S]) = i.$$
In particular, the Euler class $[e] \in H^2(\diffeo^+(\R^2);\Z)$ is
{\em unbounded}.
\end{Euler_class_thm}

This answers a question of Bestvina.

Using this result, we are able to construct examples of finitely-generated
torsion-free groups of orientation-preserving homeomorphisms of $\R^2$
which are not circularly-orderable,
thereby answering a question of Farb.

It might seem from this theorem that there are no homological constraints
on
group actions on $\R^2$, but in fact for $C^1$ actions, we show the
following:

\begin{Z+Z_thm}
Let $\rho\co\Z \oplus \Z \to \homeo^+(\R^2)$ be a $C^1$ action.
Then the Euler class $\rho^*([e]) \in H^2(\Z\oplus \Z;\Z)$ is zero.
\end{Z+Z_thm}

It should be emphasized that this is
{\em not} a purely local theorem, but uses in an essential way Brouwer's
famous theorem on fixed-point-free orientation-preserving homeomorphisms
of $\R^2$.

Together with an example of Bestvina,
theorems C and D give a complete homological classification
of (orientation-preserving) actions of (oriented) surface groups on
$\R^2$ in every degree of smoothness.

\subsection{Acknowledgements}

I would like to thank Mladen Bestvina, Nathan Dunfield, Bob Edwards,
Benson Farb,
John Franks, \'Etienne Ghys, Michael Handel, Dale Rolfsen, Fr\'ed\'eric
le Roux,
Takashi Tsuboi, Amie Wilkinson and the anonymous referee
for some very useful conversations and comments.

I would especially like to single out \'Etienne Ghys for thanks, for
reading an
earlier version of this paper and providing me with copious comments,
observations,
references, and counterexamples to some naive conjectures.

While writing this paper, I received partial support from the Sloan
foundation, and from NSF grant DMS-0405491.

\section{Left-orderable groups and circular groups}

In this section we define left-orderable and circularly-orderable
groups,
and present some of their elementary properties. None of the material in
this section is new, but perhaps the exposition will be useful to
the reader.
Details and references can be found in \cite{Mura_Rhemtulla},
\cite{Thurston_FC2},
\cite{Ghys_87}, \cite{Maclane_homology} and \cite{Loday_cyclic}, as
well as
other papers mentioned in the text as appropriate.

\subsection{Left-invariant orders}

\begin{defn}
{\rm Let $G$ be a group. A {\em left-invariant order} on $G$ is a total
order $<$
such that, for all $\alpha,\beta,\gamma$ in $G$,
$$\alpha < \beta \text{\qua iff\qua} \gamma \alpha < \gamma \beta.$$
A group which admits a left-invariant order is said to be {\em
left-orderable}.}
\end{defn}

We may sometimes abbreviate ``left-orderable'' to LO.
Note that a left-order\-able group may admit many distinct left-invariant
orders.
For instance, the group $\Z$ admits exactly two left-invariant orders.

The following lemma gives a characterization of left-orderable groups:

\begin{lem}\label{partition}
A group $G$ admits a left-invariant order iff there is a disjoint
partition
of $G = P \cup N \cup \id$ such that $P \cdot P \subset P$ and $P^{-1}
= N$.
\end{lem}
\begin{proof}
If $G$ admits a left-invariant order, set $P = \lbrace g \in G : g >
\id \rbrace$.
Conversely, given a partition of $G$ into $P,N,\id$ with the properties
above,
we can define a left-invariant order by setting $h < g$ iff $h^{-1}g
\in P$.
\end{proof}

Notice that Lemma~\ref{partition} implies that any nontrivial LO group is
{\em infinite}, and {\em torsion-free}. Notice also that any partition of
$G$ as in Lemma~\ref{partition} satisfies $N \cdot N \subset N$. For such
a partition, we sometimes refer to $P$ and $N$ as the {\em positive} and
{\em negative cone} of $G$ respectively.

LO is a {\em local} property. That is to say, it depends only on the
{\em finitely-generated} subgroups of $G$.
We make this precise in the next two lemmas.
First we show that if a group fails to be left-orderable, this fact can
be verified by examining a {\em finite} subset of the multiplication
table for
the group, and applying the criterion of Lemma~\ref{partition}.

\begin{lem}\label{finite_bad_set_LO}
A group $G$ is not left-orderable iff there is some {\em finite} symmetric
subset $S=S^{-1}$ of $G$ with the property that for every disjoint
partition
$S \backslash \id = P_S \cup N_S$, one of the following two properties
holds:
\begin{enumerate}
\item{$P_S \cap {P_S}^{-1} \ne \emptyset$ or $N_S \cap {N_S}^{-1}
\ne \emptyset$}
\item{$(P_S \cdot P_S) \cap N_S \ne \emptyset$ or
$(N_S \cdot N_S) \cap P_S \ne \emptyset$}
\end{enumerate}
\end{lem}
\begin{proof}
It is clear that the existence of such a subset contradicts
Lemma~\ref{partition}. So it suffices to show the converse.

The set of partitions of $G\backslash \id$ into disjoint sets $P,N$
is just
$2^{G\backslash \id}$, which is compact with the product topology by
Tychonoff's theorem.
By abuse of notation, if $\pi \in 2^{G\backslash \id}$ and $g \in
G\backslash \id$, we
write $\pi(g) = P$ or $\pi(g) = N$ depending on whether the element $g$
is put into the set $P$ or $N$ under the partition corresponding to $\pi$.

For every element $\alpha \in G\backslash \id$, define $A_\alpha$ to be
the {\em open}
subset of $2^{G\backslash \id}$ for which $\pi(\alpha) =
\pi(\alpha^{-1})$.
For every pair of elements $\alpha,\beta \in G\backslash \id$ with
$\alpha \ne \beta^{-1}$, define $B_{\alpha,\beta}$
to be the {\em open} subset of $2^{G\backslash \id}$ for which
$\pi(\alpha) = \pi(\beta)$
but $\pi(\alpha) \ne \pi(\alpha\beta)$.

Now, if $G$ is not LO, then by Lemma~\ref{partition}, every partition
$\pi \in 2^{G\backslash \id}$ is contained in some $A_\alpha$ or
$B_{\alpha,\beta}$.
That is, the sets $A_\alpha,B_{\alpha,\beta}$ are an open cover of
$2^{G\backslash \id}$.
By compactness, there is some {\em finite} subcover. Let $S$ denote the
set of
indices of the sets $A_\alpha,B_{\alpha,\beta}$ appearing in this finite
subcover,
together with their inverses. Then $S$ satisfies the statement of
the lemma.
\end{proof}

\begin{rmk}\label{positive_products}
{\rm An equivalent statement of this lemma is that for a group $G$ which is
not LO, there is a finite subset
$S =\lbrace g_1, \cdots, g_n \rbrace \subset G\backslash \id$ with
$S \cap S^{-1} = \emptyset$
such that for all choices of signs $e_i \in \pm 1$, the semigroup
generated
by the $g_i^{e_i}$ contains $\id$.

To see this, observe that a choice of sign $e_i \in \pm 1$ amounts to
a choice
of partition of $S \cup S^{-1}$ into $P_S$ and $N_S$. Then if $G$ is
not LO,
the semigroup of positive products of the $P_S$ must intersect the
semigroup of
positive products of the $N_S$; that is, $p = n$ for $p$ in the semigroup
generated
by $P_S$ and $n$ in the semigroup generated by $N_S$. But this implies
$n^{-1}$ is
in the semigroup generated by $P_S$, and therefore so too is the product
$n^{-1} p = \id$.}
\end{rmk}

\begin{rmk}
{\rm Given a finite symmetric subset $S$ of $G$ and a multiplication table
for $G$,
one can check by hand whether the set $S$ satisfies the hypotheses of
Lemma~\ref{finite_bad_set_LO}. It follows that if $G$ is a group
for which there is an algorithm to solve the word problem, then if
$G$ is not left-orderable, one can certify that $G$ is not left-orderable
by a
finite combinatorial certificate.}
\end{rmk}

The next lemma follows easily from Lemma~\ref{finite_bad_set_LO}:

\begin{lem}\label{local_LO}
A group $G$ is left-orderable iff every finitely-generated subgroup is
left-orderable.
\end{lem}
\begin{proof}
We use the $A,B$ notation from Lemma~\ref{finite_bad_set_LO}.

Observe that a left-ordering on $G$ restricts
to a left-ordering on any finitely-generated
subgroup $H < G$.

Conversely, suppose $G$ is not left-orderable. By
Lemma~\ref{finite_bad_set_LO}
we can find a finite set $S$ satisfying the hypotheses of that lemma. Let
$H$ be the group generated by $S$. Then Lemma~\ref{finite_bad_set_LO}
implies
that $H$ is not left-orderable.
\end{proof}

\begin{rmk}
{\rm To see this in more topological terms: observe that there is a restriction
map
$$\text{res}\co 2^{G\backslash \id} \to 2^{H \backslash \id}$$
which is surjective, and continuous with respect to the product
topologies.
It follows that the  union of the sets
$\text{res}(A_\alpha),\text{res}(B_{\alpha,\beta})$ with
$\alpha,\beta \in S$ is an open cover of $2^{H \backslash \id}$,
and therefore
$H$ is not left-orderable.}
\end{rmk}

We now study homomorphisms between LO groups.

\begin{defn}
{\rm Let $S$ and $T$ be totally-ordered sets. A map $\phi\co S \to T$ is {\em
monotone}
if for every pair $s_1,s_2 \in S$ with $s_1 > s_2$, either $\phi(s_1)
> \phi(s_2)$
or $\phi(s_1) = \phi(s_2)$.

Let $G$ and $H$ be left-orderable groups, and choose a left-invariant
order on each of
them. A homomorphism $\phi\co G \to H$ is {\em monotone} if it is monotone
as a map
or totally-ordered sets.}
\end{defn}

LO behaves well under short exact sequences:

\begin{lem}\label{ses_LO}
Suppose $K,H$ are left-orderable groups, and suppose we have a short
exact sequence
$$0 \longrightarrow K \longrightarrow G \longrightarrow H \longrightarrow 0.$$
Then for every left-invariant order on $K$ and $H$,
the group $G$ admits a left-invariant order compatible with that of $K$,
such that the
surjective homomorphism to $H$ is monotone.
\end{lem}
\begin{proof}
Let $\phi\co G \to H$ be the homomorphism implicit in the short exact
sequence.
The order on $G$ is uniquely determined by the properties that it is
required to
satisfy:
\begin{enumerate}
\item{If $\phi(g_1) \ne \phi(g_2)$ then $g_1 > g_2$ in $G$
iff $\phi(g_1) > \phi(g_2)$ in $H$}
\item{If $\phi(g_1) = \phi(g_2)$ then $g_2^{-1}g_1 \in K$, so $g_1 >
g_2$ in $G$
iff $g_2^{-1}g_1 > \id$ in $K$}
\end{enumerate}
This defines a total order on $G$ and is left-invariant, as required.
\end{proof}

\begin{defn}
{\rm A group $G$ is {\em locally LO--surjective} if every finitely-gener\-ated
subgroup $H$ admits a surjective homomorphism $\phi_H\co H \to L_H$
to an infinite LO group $L_H$.

A group $G$ is {\em locally indicable} if every finitely-generated
subgroup $H$
admits a surjective homomorphism to $\Z$. In particular, a locally
indicable group
is locally LO--surjective, though the converse is not true.}
\end{defn}

The following theorem is proved in \cite{Burns_Hale}. We give a sketch
of a proof.

\begin{thm}[Burns--Hale]\label{locally_surjective_LO}
Suppose $G$ is locally LO--surjective. Then $G$ is LO.
\end{thm}
\begin{proof}
Suppose $G$ is locally LO--surjective but not LO. Then by
Remark~\ref{positive_products},
there is some finite subset $\lbrace g_1,\dots, g_n \rbrace \subset G
\backslash \id$
such that, for all choices of signs $e_i \in \pm 1$, the semigroup of
positive products of the elements $g_i^{e_i}$ contains $\id$.
Choose a set of such $g_i$ such that $n$ is smallest possible (obviously,
$n \ge 2$).
Let $G' = \langle g_1,\dots,g_n \rangle$. Then $G'$ is finitely-generated.
Since $G$ is locally LO--surjective,
$G'$ admits a surjective homomorphism to an infinite LO group
$$\varphi\co G' \to H$$
with kernel $K$. By the defining property of the $\lbrace g_i \rbrace$,
at least one $g_i$ is in $K$ since otherwise there exist choices of
signs $e_i \in \pm 1$ such that $\varphi(g_i^{e_i})$ is in the positive
cone
of $H$, and therefore the same is true for the semigroup of positive
products
of such elements. But this would imply that the semigroup of positive
products of
the $g_i^{e_i}$ does not contain $\id$ in $G'$, contrary to assumption.
Furthermore, since $H$
is nontrivial and $\varphi$ is surjective, at least one $g_j$ is not
in $K$.

Reorder the indices of the $g_i$ so that $g_1,\dots,g_k \notin K$ and
$g_{k+1},\dots,g_n \in K$. Let $P(H)$ denote the positive elements of $H$.
Since the $g_i$ with $i \le k$ are
not in $K$, it follows that
there are choices $\delta_1,\dots,\delta_k \in \pm 1$ such that
$\varphi(g_i^{\delta_i}) \in P(H)$. Moreover, since $n$ was chosen
to be minimal, there exist choices $\delta_{k+1}, \dots, \delta_n \in
\pm 1$ such
that no positive product of elements of $g_{k+1}^{\delta_{k+1}},\dots,
g_n^{\delta_n}$
is equal to $\id$.

On the other hand, by the definition of $g_i$, there are positive integers
$n_i$ such that
$$\id = g_{i(1)}^{n_1\delta_{i(1)}} \cdots g_{i(s)}^{n_s\delta_{i(s)}}$$
where each $i(j)$ is between $1$ and $n$. By hypothesis, $i(j) \le k$
for at least
one $j$. But this implies that the image of the right hand side of
this equation
under $\varphi$ is in $P(H)$, which is a contradiction.
\end{proof}

Theorem~\ref{locally_surjective_LO} has the corollary that a locally
indicable group is LO. It is this corollary that will be most useful
to us.

\subsection{Circular orders}

The approach we take in this section is modelled on \cite{Thurston_FC2},
although
an essentially equivalent approach is found in \cite{Ghys_87}.

We first define a circular ordering on a set. Suppose $p$ is a point in an
oriented circle $S^1$. Then $S^1 \backslash p$ is homeomorphic to $\R$,
and the
orientation on $\R$ defines a natural total order on $S^1 \backslash p$.
In general, a circular order on a set $S$ is defined by a choice of total
ordering on each subset of the form $S \backslash p$, subject to certain
compatibility
conditions which we formalize below.

\begin{defn}
{\rm Let $S$ be a set. A {\em circular ordering} on a set $S$ with at least $4$
elements is a choice of
total ordering on $S \backslash p$ for every $p \in S$, such that if
$<_p$ is the total ordering defined by $p$, and $p,q \in S$ are two
distinct
elements, the total orderings $<_p,<_q$ differ by a {\em cut} on their
common
domain of definition. That is, for any $x,y$ distinct from $p,q$, the
order of
$x$ and $y$ with respect to $<_p$ and $<_q$ is the same unless
$x <_p q <_p y$, in which case we have $y <_q p <_q x$.
We also say that the order $<_q$ on $S\backslash \lbrace p, q\rbrace$
is obtained from the
order $<_p$ on $S \backslash p$ by {\em cutting at $q$}.

If $S$ has exactly $3$ elements $S = \lbrace x,y,z \rbrace$,
we must add the condition that $y <_x z$ iff $z <_y x$. Note that this
condition
is implied by the condition in the previous paragraph if $S$ has at
least $4$ elements.}
\end{defn}

To understand the motivation for the terminology, consider the operation
of
{\em cutting} a deck of cards.

\begin{exa}
{\rm The oriented circle $S^1$ is circularly-ordered, where for any $p$, the
ordering $<_p$ is just the ordering on $S^1 \backslash p \cong \R$
induced by
the orientation on $\R$.}
\end{exa}

\begin{defn}
{\rm A set with three elements $x,y,z$ admits exactly two circular orders,
depending
on whether $y <_x z$ or $z <_x y$. In the first case, we say the triple
$(x,y,z)$ is {\em positively-ordered} and in the second case, we say it is
{\em negatively-ordered}.}
\end{defn}

We also refer to a positively-ordered triple of points as {\em
anticlockwise}
and a negatively-ordered triple as {\em clockwise}, by analogy with the
standard circular order on triples of points in the positively oriented
circle.

A circular ordering on a set $S$ induces a circular ordering on any subset
$T \subset S$. If $T_\alpha$ is a family of subsets of $S$ which are
all circularly-ordered, we say the circular orderings on the $T_\alpha$
are {\em compatible} if
they are simultaneously induced by some circular ordering on $S$.

It is clear that a circular ordering on a set $S$ is determined by the
family of circular orderings on all triples of elements in
$S$. Conversely,
the following lemma characterizes those families of circular orderings
on triples
of elements which arise from a circular ordering on all of $S$:

\begin{lem}\label{compatible_on_quadruples}
Suppose $S$ is a set. A circular ordering on all triples of distinct
elements on $S$
is compatible iff for every subset $Q \subset S$ with four elements, the
circular ordering on triples of distinct elements of $Q$ is compatible. In
this
case, these circular orderings are {\em uniquely} compatible, and
determine a circular
ordering on $S$.
\end{lem}
\begin{proof}
A circular ordering on triples in $S$ defines, for any $p \in S$,
a binary relation
$<_p$ on $S\backslash p$ by $x <_p y$ iff the triple $(p,x,y)$ is
positively-ordered.
To see that this binary relation defines a total ordering on $S\backslash
p$,
we must check transitivity of $<_p$. But this follows from compatibility
of
the circular ordering on quadruples $Q$. It is straightforward to
check that the total orders $<_p$ and  $<_q$
defined in this way differ by a cut for distinct $p,q$.
\end{proof}

\begin{defn}
{\rm Let $C_1,C_2$ be circularly-ordered sets. A map $\phi\co C_1 \to C_2$ is
{\em monotone} if for each $c \in C_2$ and each $d \in \phi^{-1}(c)$,
the restriction map between totally-ordered sets
$$\phi\co (C_1 \backslash \phi^{-1}(c), <_d) \to (C_2 \backslash c, <_c)$$
is monotone.}
\end{defn}

There is a natural topology on a circularly-ordered set for which monotone
maps are continuous.

\begin{defn}
{\rm Let $O, <$ be a totally-ordered set. The {\em order topology} on $O$
is the
topology generated by open sets of the form $\lbrace x | x > p \rbrace$
and $\lbrace x | x < p \rbrace$ for all $p \in O$.
Let $S$ be a circularly-ordered set. The {\em order topology} on $S$ is
the topology generated on each $S\backslash p$ by the (usual) order
topology
on the totally-ordered set $S\backslash p, <_p$.}
\end{defn}

We now turn to the analogue of left-ordered groups for circular orderings.

\begin{defn}
{\rm A group $G$ is {\em left circularly-ordered} if it admits a circular
order as
a set which is preserved by the action of $G$ on itself on the left.
A group is {\em left circularly-orderable} if it can be left
circularly-ordered.}
\end{defn}

We usually abbreviate this by saying that a group is {\em
circularly-orderable}
if it admits a {\em circular order}.

\begin{exa}
{\rm A left-orderable group $G,<$ is circularly-orderable as follows: for each
element $g \in G$, the total order $<_g$ on $G\backslash g$ is obtained
from
the total order $<$ by cutting at $g$.}
\end{exa}

\begin{defn}
{\rm The group of orientation-preserving homeomorph\-isms of $\R$ is denoted
$\homeo^+(\R)$. The group of orientation-preserving homeomorph\-isms
of the circle is denoted $\homeo^+(S^1)$.}
\end{defn}

An action of $G$ on $\R$ or the circle by orientation-preserving
homeomorph\-isms is the
same thing as a representation in $\homeo^+(\R)$ or $\homeo^+(S^1)$.
We will see that for countable groups $G$, being LO is the same as
admitting a faithful representation in $\homeo^+(\R)$, and
CO is the same as admitting a faithful representation in $\homeo^+(S^1)$.
First we give one direction of the implication.

\begin{lem}\label{CO_group_acts}
If $G$ is countable and admits a left-invariant circular order, then
$G$ admits a faithful representation in $\homeo^+(S^1)$.
\end{lem}
\begin{proof}
Let $g_i$ be a countable enumeration of the elements of $G$.
We define an embedding $e\co G \to S^1$ as follows.
The first two elements $g_1,g_2$ map to arbitrary distinct
points in $S^1$. Thereafter, we use the following inductive
procedure to uniquely extend $e$ to each $g_n$.

Firstly, for every $n>2$, the map
$$e\co \bigcup_{i \le n} g_i \longrightarrow \bigcup_{i \le n} e(g_i)$$
should be injective and circular-order-preserving, where the $e(g_i)$ are
circul\-arly-ordered by the natural circular ordering on $S^1$.
Secondly, for every $n>2$, the element $e(g_n)$ should be taken to the
midpoint of the unique interval complementary to $\bigcup_{i < n}
e(g_i)$ compatible
with the first condition. This defines $e(g_n)$ uniquely, once $e(g_i)$
has
been defined for all $i<n$.

It is easy to see that the left action of $G$ on itself
extends uniquely to a continuous order preserving homeomorphism
of the closure $\overline{e(G)}$ to itself. The complementary intervals
$I_i$ to
$\overline{e(G)}$ are permuted by the action of $G$; we choose an
identification
$\varphi_i\co I_i \to I$ of each interval with $I$, and extend the
action of
$G$ so that if $g(I_i) = I_j$ then the action of $g$ on $I_i$ is equal to
$$g|_{I_i} = \varphi_j^{-1}\varphi_i.$$
This defines a faithful representation of $G$ in $\homeo^+(S^1)$,
as claimed.
\end{proof}

\begin{rmk}
{\rm Note that basically the same argument shows that a left-order\-able
countable
group is isomorphic to a subgroup of $\homeo^+(\R)$.
Notice further that this construction has an important property: if
$G$ is a countable left- or circularly-ordered group, then $G$
is circular
or acts on $\R$ in such a way that {\em some point has trivial
stabilizer}.
In particular, any point in the image of $e$ has trivial stabilizer.}
\end{rmk}

Short exact sequences intertwine circularity and left-orderability:
\begin{lem}\label{ses_CO}
Suppose
$$0 \longrightarrow K \longrightarrow G \longrightarrow H \longrightarrow 0$$
is a short exact sequence, where $K$ is
left-ordered and $H$ is circularly-ordered. Then $G$ can be
circularly-ordered
in such a way that the inclusion of $K$ into $G$ respects the order on
$G\backslash g$ for any $g$ not in $K$, and the map from
$G$ to $H$ is monotone.
\end{lem}
\begin{proof}
Let $\phi\co G \to H$ be the homomorphism in the short exact sequence.
Let $g_1,g_2,g_3$ be three distinct elements of $G$. We define the
circular order as
follows:
\begin{enumerate}
\item{If $\phi(g_1),\phi(g_2),\phi(g_3)$ are distinct, circularly-order
them by the
circular order on their image in $H$}
\item{If $\phi(g_1)=\phi(g_2)$ but these are distinct from $\phi(g_3)$,
then
$g_2^{-1}g_1 \in K$. If $g_2^{-1}g_1 < \id$ then $g_1,g_2,g_3$ is
positively-ordered, otherwise it is nega\-tively-ordered}
\item{If $\phi(g_1)=\phi(g_2)=\phi(g_3)$ then
$g_3^{-1}g_1,g_3^{-1}g_2,\id$ are all in $K$,
and therefore inherit a total ordering. The three corresponding elements
of $G$ in the
same total order are negatively-ordered}
\end{enumerate}
One can check that this defines a left-invariant circular order on $G$.
\end{proof}

Here our convention has been that the orientation-preserving
inclusion of $\R$ into $S^1 \backslash p$ is order-preserving.

We will show that for countable groups, being LO or CO is equivalent
to admitting
a faithful representation in $\homeo^+(\R)$ or $\homeo^+(S^1)$
respectively.
But first we must describe an operation due to Denjoy \cite{Denjoy_blowup}
of {\em blowing up} or {\em Denjoying} an action.

\begin{construct}[Denjoy]\label{blow_up}
{\rm Let $\rho\co G \to \homeo^+(S^1)$ be an action of a countable group on
$S^1$. For convenience, normalize $S^1$ to have length $1$.
Let $p \in S^1$ be some point. Let $O$ denote the countable orbit of
$p$ under $G$, and let $\phi\co O \to \R^+$ assign a positive real
number to
each $o \in O$ such that $\sum_{o \in O} \phi(o) = 1$.
Choose some point $q$ not in $O$, and define $\tau\co [0,1] \to S^1$ to be
an orientation-preserving parametrization by length,
which takes the two endpoints to $q$. Define $\sigma\co [0,1] \to
[0,2]$ by
$$\sigma(t) = t + \sum_{o \in O~:~\tau^{-1}(o) \le t} \phi(o).$$
Then $\sigma$ is discontinuous on $\tau^{-1}(O)$, and its graph can be
completed to
a continuous image of $I$ in $[0,1] \times [0,2]$
by adding a vertical segment of length $\phi(o)$ at each
point $\tau^{-1}(o)$ where $o \in O$. Identify opposite sides of $[0,1]
\times [0,2]$
to get a torus, in which the closure of the graph of $\sigma$ closes up
to become
a $(1,1)$ curve which, by abuse of notation, we also refer to as $\sigma$.
Notice that projection $\pi_h$ onto the horizontal factor defines
a monotone
map from $\sigma$ to $S^1$.

Then the action of $G$ on $S^1$ extends in an obvious way
to an action on this torus which leaves the
$(1,1)$ curve invariant, and also preserves the foliations of the torus
by horizontal
and vertical curves. Up to conjugacy in $\homeo^+(\sigma)$, the action
of $G$ on
$\sigma$ is well-defined, and is called the {\em blown-up action
at $p$}.
The pushforward of this blown-up action under $(\pi_h)_*$ recovers
the original
action of $G$ on $S^1$; that is, the two actions are related by a
degree one
{\em monotone map}, and are said to be {\em semi-conjugate}. The
equivalence
relation that this generates is called {\em monotone equivalence}.}
\end{construct}

With this construction available to us, we demonstrate the equivalence
of CO
with admitting a faithful representation in $\homeo^+(S^1)$.

\begin{thm}\label{action_is_order}
A countable group $G$ is left- or circularly-ordered
iff $G$ admits a
faithful homomorphism to $\homeo^+(\R)$ or $\homeo^+(S^1)$ respectively.
Moreover, the action on $\R$ or $S^1$ can be chosen so that some point
has a
trivial stabilizer.
\end{thm}
\begin{proof}
In Lemma~\ref{CO_group_acts} we have already showed
how a left or circular order gives rise to a faithful
action on $\R$ or $S^1$. So it remains to prove the converse.

Let $\phi\co G \to \homeo^+(\R)$ be faithful. Let $p_i$ be some sequence
of points such
that the intersection of the stabilizers of the $p_i$ is the
identity. Some such sequence
$p_i$ exists, since $G$ is countable, and any nontrivial element acts
nontrivially
at some point. Then each $p_i$ determines a (degenerate) left-invariant
order on $G$,
by setting $g>_i h$ if $g(p_i)> h(p_i)$, and $g=_i h$ if $g(p_i) =
h(p_i)$. Then we define
$g>h$ if $g >_i h$ for some $i$, and $g =_j h$ for all $j<i$.

The definition of a circular order is similar: pick some point $p \in
S^1$, and
suppose that the stabilizer $\stab(p)$ is nontrivial. Then $\stab(p)$ acts
faithfully on $S^1 - p = \R$, so by the argument above, $\stab(p)$
is left-orderable
and acts on $\R$. In fact, we know $\stab(p)$ acts on $\R$ in such a
way that
some point has trivial stabilizer. Let $\varphi\co \stab(p) \to
\homeo^+(\R)$
be such
a representation. We construct a {\em new} representation
$\phi'\co G \to \homeo^+(S^1)$ from $\phi$ by {\em blowing up $p$}
as in Construction~\ref{blow_up}. The representation $\phi'$ is monotone
equivalent
to $\phi$; that is, there is a monotone map $\pi\co S^1 \to S^1$
satisfying
$$\pi_*\phi' = \phi.$$
Let $C \subset S^1$ be the set where the monotone map $\pi$ is not
locally constant.
We will modify the action of $G$ on $S^1\backslash C$ as follows. Note
that
$G$ acts on $C$ by the pullback under $\pi$
of the action on $S^1$ by $\phi$. We extend this action to $S^1\backslash
C$ to define
$\phi''$. Let $I$ be the open interval obtained by blowing up $p$. We
identify $I$ with
$\R$, and then let $\stab(p)$ act on $I$ by the pullback of $\varphi$
under this
identification. Each other component $I_i$ in
$S^1\backslash C$ is of the form $g(I)$ for some $g \in G$. Choose such
a $g_i$ for each
$I_i$, and pick an arbitrary (orientation preserving) identification
$\varphi_i\co I \to I_i$,
and define $\phi''(g_i)|_I = \varphi_i$. Now, for any $g \in G$, define
$g|_{I_i}$ as
follows: suppose $g(I_i) = I_j$. Then $g_j^{-1}gg_i \in \stab(p)$,
so define
$$\phi''(g)|_{I_i} = \varphi_j \varphi(g_j^{-1}gg_i) \varphi_i^{-1}\co I_i
\to I_j.$$
It is clear that this defines a faithful representation $\phi''\co G
\to \homeo^+(S^1)$,
monotone equivalent to $\phi$, with the property that some point $q \in
S^1$ has trivial
stabilizer.

Now define a circular order on distinct
triples $g_1,g_2,g_3$ by restricting the circular order on $S^1$ to the
triple $g_1(q),g_2(q),g_3(q)$.
\end{proof}

Notice that in this theorem, in order to recover a left or circular
order on $G$ from a
faithful action, all we used about $\R$ and $S^1$ was that they were
ordered
and circularly-ordered sets respectively.

With this theorem, and our lemmas on short exact sequences,
we can deduce the existence of left or circular orders on countable
groups from
the existence of actions on ordered or circularly-ordered sets, with
left-orderable kernel.

\begin{thm}\label{action_kernel}
Suppose a countable group $G$ admits an action by order preserving maps
on a totally-ordered or circularly-ordered set $S$ in such a way that
the kernel $K$ is left-orderable. Then $G$ admits a faithful, order
preserving action on $\R$ or $S^1$, respectively.
\end{thm}
\begin{proof}
We discuss the case that $S$ is circularly-ordered, since this is slightly
more complicated. Since $G$ is countable, it suffices to look at an orbit
of the action, which will
also be countable. By abuse of this notation, we also denote the orbit by
$S$. As in Lemma~\ref{CO_group_acts}, the set $S$ with its order topology
is naturally order-isomorphic to a subset of $S^1$. Let $\overline{S}$
denote
the closure of $S$ under this identification. Then the action of $G$ on
$S$ extends to an orientation-preserving action on $S^1$, by permuting
the
complementary intervals to $\overline{S}$. It follows that the image of
$G$ in $\homeo^+(S^1)$ is CO, with kernel $K$. By Lemma~\ref{ses_CO},
$G$ is CO.
By Theorem~\ref{action_is_order}, the proof follows.

The construction for $S$ totally-ordered is similar.
\end{proof}

\subsection{Homological characterization of circular groups}

Circular orders on groups $G$ can be characterized homologically. There
are
at least two different ways of doing this, due to Thurston and Ghys
respectively,
which reflect two different ways of presenting the theory of group
cohomology.

First, we recall the definition of group cohomology. For details, we
refer to
\cite{Maclane_homology} or \cite{Loday_cyclic}.

\begin{defn}
{\rm Let $G$ be a group.
The {\em homogeneous chain complex} of $G$ is a complex $C_*(G)_h$ where
$C_n(G)_h$ is the free abelian group generated by equivalence classes of
$(n+1)$--tuples $(g_0:g_1:\dots:g_n)$, where two such tuples are
equivalent if they are in the same coset of the left diagonal action
of $G$
on the coordinates. That is,
$$(g_0:g_1:\dots:g_n) \sim (gg_0:gg_1:\dots:gg_n).$$
The boundary operator in homogeneous coordinates is very simple,
defined by the
formula
$$\partial (g_0:\dots:g_n) = \sum_{i=0}^n (-1)^i (g_0:
\dots:\hat{g_i}:\dots:g_n).$$
The {\em inhomogeneous chain complex} of $G$ is a complex $C_*(G)_i$
where $C_n(G)_i$ is the free abelian group generated by $n$--tuples
$(f_1,\dots,f_n)$. The boundary operator in inhomogeneous coordinates is
more complicated, defined by the formula
\begin{align*}
\partial (f_1,\dots,f_n) & = (f_2,\dots,f_n) + \sum_{i=1}^{n-1} (-1)^i
(f_1,\dots,f_i f_{i+1},\dots,f_n) \\
& + (-1)^n (f_1,\dots,f_{n-1}).
\end{align*}
The relation between the two coordinates comes from the following
bijection
of generators
$$(g_0:g_1:\dots:g_n) \longrightarrow
(g_0^{-1}g_1,g_1^{-1}g_2,\dots,g_{n-1}^{-1}g_n)$$
which correctly transforms one definition of $\partial$ to the other.
It follows that the two chain complexes are canonically isomorphic, and
therefore by abuse of notation we denote either by $C_*(G)$, and
write an element either in homogeneous or inhomogeneous coordinates
as convenient.

Let $R$ be a commutative ring. The homology of the complex $C_*(G)\otimes
R$
is denoted $H_*(G;R)$, and the homology of the
adjoint complex $\hom(C_*(G),R)$ is denoted $H^*(G;R)$. If $R = \Z$,
we abbreviate
these groups to $H_*(G)$ and $H^*(G)$ respectively. If $G$ is a
topological
group, and we want to stress that this is the abstract group (co)homology,
we
denote these groups by $H_*(G^\delta)$ and $H^*(G^\delta)$ respectively
($\delta$ denotes the discrete topology).}
\end{defn}

We give a geometrical interpretation of this complex.
The simplicial realization of the complex $C_*(G)$ is
a model for the classifying space $BG$, where $G$ has the discrete
topology.
If $G$ is torsion-free, an equivalent model for $EG$ is the
complete simplex on the elements of $G$.
In this case, since $G$ is torsion-free,
it acts freely and properly discontinuously on this simplex, with
quotient $BG$. If we label vertices of $EG$ tautologically by elements
of $G$,
the labels on each simplex give
homogeneous coordinates on the quotient. If we label edges of $EG$
by the
difference of the labels on the vertices at the ends, then the labels are
well-defined on the quotient; the labels on the $n$ edges between
consecutive
vertices of an $n$--simplex,
with respect to a total order of the vertices, give inhomogeneous
coordinates.

In particular, the cohomology $H^*(G)$ is just the cohomology of the
$K(G,1)$, that is, of the unique (up to homotopy) aspherical space with
fundamental group isomorphic to $G$. If $G$ is not torsion-free,
this equality
of groups is nevertheless true.

The cohomology of the group $\homeo^+(S^1)$ is known by a general theorem
of Mather
and Thurston (see \cite{Thurston_diffeo} or \cite{Tsuboi} for details
and more references):

\begin{thm}[Mather, Thurston]\label{cohomology_of_homeos}
For any manifold $M$, there is an isomorphism of cohomology rings
$$H^*(\homeo(M)^\delta;\Z) \cong H^*(\bhomeo(M);\Z)$$
where $\bhomeo(M)$ denotes the classifying space of the {\em topological
group}
of homeomorphisms of $M$, and the left hand side denotes the {\em group
cohomology}
of the {\em abstract} group of homeomorphisms of $M$.
\end{thm}

For $M = S^1$ or $\R^2$, the topological group $\homeo^+(M)$ is homotopy
equivalent to a circle. For $S^1$, this is trivial. For $\R^2$, we
observe that $\homeo^+(\R^2)$ is the stabilizer of a point in
$\homeo^+(S^2)$, and then apply a theorem of Kneser \cite{Kneser} about
the homotopy type of $\homeo^+(S^2)$.  It follows that $\bhomeo^+(M)$
in either case is homotopic to $\CP^\infty$, and therefore there is an
isomorphism of rings $$H^*(\homeo^+(\R^2);\Z) \cong H^*(\homeo^+(S^1);\Z)
\cong \Z[e]$$ where $[e]$ is a free generator in degree $2$ called the
{\em Euler class}.

An algebraic characterization of the Euler class can be given.

\begin{defn}
{\rm For any group $G$ with $H^1(G;\Z) = 1$, there is a {\em universal
central extension}
$$0 \longrightarrow A \longrightarrow \hat{G} \longrightarrow G
  \longrightarrow 0$$
where $A$ is abelian, with the property that for any other central
extension
$$0 \longrightarrow B \longrightarrow G' \longrightarrow G \longrightarrow 0$$
there is a unique homomorphism from $\hat{G} \to G'$, extending uniquely
to a
morphism of short exact sequences.}
\end{defn}

A non-split central extension $G'$ can be characterized as the universal
central extension
of $G$ iff $G$ is perfect (i.e. $H^1(G;\Z) = 1$) and every central
extension of
$G'$ splits. See Milnor \cite{Milnor_Ktheory} for more details.

For $G = \homeo^+(S^1)$, the universal central extension is denoted
%${\thomeo^+(S^1)}$,
$\thomeo^+(S^1)$,
and can be identified with the preimage of $\homeo^+(S^1)$ in
$\homeo^+(\R)$ under
the covering map $\R \to S^1$. The center of ${\thomeo^+(S^1)}$
is $\Z$, and
the class of this $\Z$ extension is called the Euler class. By the
universal
property of this extension, one sees that this class is the generator of
$H^2(\homeo^+(S^1);\Z)$. This can be summarized by a short exact sequence
$$0 \longrightarrow \Z \longrightarrow {\thomeo^+(S^1)} \longrightarrow
  \homeo^+(S^1) \longrightarrow 0.$$
The following construction is found in \cite{Thurston_FC2}. An equivalent
construction is given in \cite{Jekel}.

\begin{construct}[Thurston]\label{Thurston_cocycle}
{\rm Let $G$ be a countable CO group, and let $\rho\co G \to \homeo^+(S^1)$ be
constructed as in Theorem~\ref{action_is_order} so that the point $p$
has trivial stabilizer. For each triple $g_0,g_1,g_2 \in G$ of
distinct elements, define the cocycle
$$c(g_0:g_1:g_2) = \left\{ \begin{array}{rl}
  1 & \text{ if } (g_0(p),g_1(p),g_2(p)) \text{ is positively oriented} \\
 -1 & \text{ otherwise.}\end{array}\right.$$
It is clear that $c$ is well-defined on the homogeneous coordinates
for $C_2(G)$.
Then extend $c$ to degenerate triples by setting it equal to $0$ if at
least two of
its coefficients are equal.

The fact that the circular order on triples of points in $S^1$ is
compatible on
quadruples is exactly the condition that the coboundary of $c$ is $0$
--- that is,
$c$ is a cocycle, and defines an element $[c] \in H^2(G;\Z)$.}
\end{construct}

The following (related) construction is found in \cite{Ghys_87}:

\begin{construct}[Ghys]\label{Ghys_cocycle}
{\rm Let $G$ be a countable CO group. Let $\rho\co G \to \homeo^+(S^1)$
be constructed
as in Theorem~\ref{action_is_order}. By abuse of notation, we identify
$G$ with
its image $\rho(G)$. Let $\hat{G}$ denote the preimage of $G$ in the
extension ${\thomeo^+(S^1)} \subset \homeo^+(\R)$.
Define a section $s\co G \to \hat{G}$ uniquely by the property that
$s(g)(0)
\in [0,1)$.
For each pair of elements $g_0,g_1 \in G$, define the cocycle
$$e(g_0,g_1) = s(g_0g_1)^{-1}s(g_0)s(g_1)(0).$$
Then one can check that $e$ is a {\em cocycle} on $C_2(G)$ in
inhomogeneous coordinates,
and defines an element $[e] \in H^2(G;\Z)$.
Moreover, $e$ takes values in $\lbrace 0,1 \rbrace$.}
\end{construct}

The following lemma can be easily verified; for a proof, we refer to
\cite{Thurston_FC2}
or \cite{Jekel}.

\begin{lem}[Ghys, Jekel, Thurston]\label{ce_euler_class}
Let $G$ be a countable circularly-ordered group.
The cocycles $e,c$ satisfy
$$2[e] = [c].$$
Moreover, the class $[e]$ is the Euler class of the circular order on $G$.
\end{lem}

Actually, the restriction to countable groups is not really necessary.
One can define the cocycles $c,e$ directly from a circular order on
an arbitrary
group $G$. This is actually done in \cite{Thurston_FC2} and
\cite{Ghys_87};
we refer the reader to those papers for the more abstract construction.

\begin{thm}\label{CO_lifts_to_LO}
Let $G$ be a circularly-ordered group with Euler class $[e] \in
H^2(G;\Z)$.
If $[e] = 0$, then $G$ is left-ordered. In any case, the central extension
of $G$
corresponding to the class $[e]$ is left-orderable.
\end{thm}
\begin{proof}
We prove the theorem for $G$ countable; the general case is proved in
\cite{Ghys_87}.

From the definition of $s$ in Construction~\ref{Ghys_cocycle} and
Lemma~\ref{ce_euler_class}, we see that $e$ is the obstruction to
finding some
(possibly different) section $G \to \hat{G}$. But $\hat{G}$ is a
subgroup of
the group $\homeo^+(\R)$. Now, every finitely-generated subgroup of
$\homeo^+(\R)$ is
left-orderable, by Theorem~\ref{action_is_order}. It follows by
Lemma~\ref{local_LO}
that the entire group $\homeo^+(\R)$ is left-orderable; in particular,
so is $\hat{G}$.
\end{proof}

\subsection{Bounded cohomology and the Milnor--Wood inequality}

Construction~\ref{Thurston_cocycle} and Construction~\ref{Ghys_cocycle}
do more
than give an explicit representative cocycle of the Euler class;
they verify
that this cocycle has a further additional property, namely that the
Euler class is a {\em bounded cocycle} on $G$.

\begin{defn}
{\rm Suppose $R=\R$ or $\Z$. Define an $L_1$ norm on $C_i(G)$ in the obvious
way by
$$\Bigg\| \sum_j s_j (g_0(j):g_1(j): \dots :g_i(j)) \Bigg\|_1 = \sum_j
|s_j|.$$
Dually, the $L_\infty$ norm is partially defined on $\hom(C_i(G);R)$,
and the
subspace consisting of homomorphisms of finite $L_\infty$ norm is denoted
$\hom_b(C_i(G);R)$.
The coboundary takes cochains of finite $L_\infty$ norms to cocycles
of finite
$L_\infty$ norm, and therefore we can take the cohomology of the
subcomplex. This
cohomology is denoted $H^*_b(G;\R)$ and is called the {\em bounded
cohomology} of $G$.
For an element $\alpha \in H^*_b(G;\R)$, the norm of $\alpha$, denoted
$\|\alpha\|_\infty$ or just $\|\alpha\|$, is the infimum of $\|c\|_\infty$
over cocycles
$c$ with $[c] = \alpha$.}
\end{defn}

When we do not make coefficients explicit, the {\em norm} of a bounded
cocycle refers
to its norm amongst representatives with $\R$ coefficients.

In this language, the famous Milnor--Wood inequality \cite{Milnor_bound},
\cite{Wood}
can be expressed as follows:

\begin{thm}[Milnor--Wood]\label{Milnor_Wood_inequality}
Let $G$ be a circularly-ordered group. Then the Euler class $[e]$ of $G$
is an
element of $H^2_b(G)$ with norm $\|[e]\| \le \frac 1 2$.
\end{thm}
\begin{proof}
Let $e$ be the cocycle constructed by Ghys. Then $\frac c 2 = e- \frac
1 2$
is homologous to $e$, and has norm $\le \frac 1 2$.
\end{proof}

We will see in Section 4 that although $\homeo^+(S^1)$ and
$\homeo^+(\R^2)$
have the same
cohomology as abstract groups, their {\em bounded cohomology} groups are
very different. This difference persists to the smooth category, as we
shall see.

\section{Planar groups with bounded orbits}

The purpose of this section is to show that every group of $C^1$
orientation-preserving homeomorphisms of $\R^2$ with a bounded orbit
is circularly-orderable. The main tools will be the Thurston stability
theorem,
and certain generalizations of the braid groups. In the course of
the proof
we also show that the mapping class group of a compact totally
disconnected
set in the plane is circularly-orderable, and the mapping class group rel.
boundary of such a set in the disk is left-orderable. This generalizes a
theorem of Dehornoy \cite{Dehorn_braids}
on orderability of the usual (finitely-generated) braid
groups.

\subsection{Prime ends}

In this section we describe some elements of the theory of prime ends.
For details of proofs and references, consult \cite{Pomm_boundary}
or \cite{Mather}.

Prime ends are a technical tool, developed in conformal analysis, to study
the boundary behaviour of conformal maps which take the unit disk in
$\C$ to
the interior $U$ of a region $K$ whose boundary is not locally
connected. They
were introduced by Carath\'eodory in \cite{Cara_prime}. If $\partial K$
{\em is}
locally connected, then the prime ends of $\partial K$ are just the
proper homotopy classes of proper rays in $U$.

\begin{defn}
{\rm Let $U$ be a bounded open subset of $\R^2$, and let $K$ denote its
closure.
Fix the notation $\partial K = K \backslash U$. Notice that this might
not be
the frontier of $K$ in the usual sense.
A {\em null chain} is a sequence of proper
arcs $(C_i,\partial C_i) \subset (K,\partial K)$ where
$\inte(C_i) \subset U$ is an embedding, but whose endpoints are not
necessarily
distinct, such that the following conditions are satisfied:
\begin{enumerate}
\item{$\overline{C_n} \cap \overline{C_{n+1}} = \emptyset$}
\item{$C_n$ separates $C_0$ from $C_{n+1}$ in $U$}
\item{The diameter of the $C_n$ converges to $0$ as $n \to \infty$.}
\end{enumerate}
A {\em prime end} of $U$ is an equivalence class of null chains, where two
null chains $\lbrace C_i \rbrace,\lbrace C_i' \rbrace$ are equivalent
if for
every sufficiently large $m$, there is an $n$ such that $C_m'$ separates
$C_0$
from $C_n$, and $C_m$ separates $C_0'$ from $C_n'$.

For such a set $U$, denote the set of prime ends of $U$ by $\P(U)$.}
\end{defn}

\begin{exa}
{\rm Let $D$ denote the open unit disk in $\C$. Then $\P(D)$
is in natural bijective correspondence with $\partial \overline{D}$.}
\end{exa}

Notice that any homeomorphism of $\R^2$ which fixes $U$ as a set will
induce an
automorphism of $\P(U)$. On the other hand, a self-homeomorphism of $U$
which does not extend continuously to $\overline{U}$ will not typically
induce an
automorphism of $\P(U)$, since for instance the image of a properly
embedded
arc $C_i$ with endpoints on $\overline{U}\backslash U$ will not
necessarily
limit to well-defined endpoints in $U$.

\begin{lem}\label{prime_ends_circular}
Let $U$ be a simply-connected, bounded, open subset of $\R^2$.
Then the set of prime ends $\P(U)$ admits a natural circular order.
\end{lem}
\begin{proof}
Let $\varphi\co U \to D$ be a uniformizing map. Then the set of prime ends
of $U$
is taken bijectively to the set of prime ends of $D$. This is not entirely
trivial; it is
contained in proposition 2.14 and theorem 2.15 in \cite{Pomm_boundary}.
In any case, this map identifies $\P(U)$ with $\partial
\overline{D}$. Thus
$\P(U)$ inherits a natural circular ordering from
$\partial \overline{D}$. If $\varphi'$ is another uniformizing map, then
$\varphi' \circ \varphi^{-1}$ is a M\"obius transformation of $\partial
\overline{D}$,
and therefore preserves the circular order.
\end{proof}

If $\varphi\co U \to V$ is a conformal map between simply connected
domains,
then let $\varphi_*\co \P(U) \to \P(V)$ denote the corresponding map
between
prime ends.
As in the proof of Lemma~\ref{prime_ends_circular}, the proof that
$\varphi_*$ is well-defined is found in \cite{Pomm_boundary}.

\begin{lem}\label{action_faithful}
Suppose $\varphi\co \R^2 \to \R^2$ fixes $U$ as a set, and
fixes the prime ends of $U$. Then $\varphi$ fixes
$\overline{U} \backslash U$ pointwise.
\end{lem}
\begin{proof}
If $f\co D \to U$ is a uniformizing map, then $f$ has a radial limit at
$\zeta \in S^1$ iff the prime end $f_*(\zeta) \in \P(U)$ is {\em
accessible};
that is, if there is a Jordan arc that lies in $U$ except for one
endpoint,
and that intersects all but finitely many crosscuts of some null-chain of
$f_*(\zeta)$. The endpoints of this Jordan arc is called an {\em
accessible point}.
It is known (see \cite{Pomm_boundary}) that the set of
accessible points is {\em dense} in $\overline{U}\backslash U$.
But an automorphism which fixes all prime ends must fix all accessible
points, and therefore must fix $\overline{U} \backslash U$ pointwise.
\end{proof}

\subsection{Groups which stabilize a point}

In this subsection we state the Thurston stability theorem, and from
this deduce
information about the group of $C^1$ homeomorphisms of $\R^2$ which
stabilizes a point.

The Thurston stability theorem is proved in
\cite{Thurston_stability}. Since
many people will only be familiar with the $1$--dimensional version of this
theorem, we indicate the idea of the proof.

\begin{thm}[Thurston stability theorem]\label{thurston_stability}
Let $p$ be a point in a smooth manifold $M^n$.
Let $G$ be a group of germs at $p$
of $C^1$ homeomorphisms of $M^n$ which fix $p$.
Let $L\co G \to \GL(n,\R)$ denote the natural homomorphism obtained
by linearizing $G$ at $p$. Let $L(G)$ denote the image of $G$, and
$K(G)$ the
kernel of $L$. Then $K(G)$ is locally indicable.
\end{thm}
\begin{proof}
The idea of the proof is as follows. Let $H < K(G)$ be a finitely
generated subgroup,
with generators $h_1,\dots h_m$. Let $p_i \to p$ be some convergent
sequence.
If we rescale the action near $p_i$ so that every $h_j$ moves points a
bounded distance,
but some $h_k(i)$ moves points distance $1$, then the rescaled actions
vary in a
precompact family. It follows that we can extract a limiting nontrivial
action,
which by construction will be an action by {\em translations}. In
particular,
$H$ is indicable, and $K(G)$ is locally indicable.

Now we make this more precise. Change coordinates so that $p$ is
at the origin. Then each generator $h_i$ can be expressed in local
coordinates as
a sum
$$h_i(x) = x + y(h_i)(x)$$
where $|y(h_i)(x)| = o(x)$ and satisfies $y(h_i)'|_0 = 0$.
For each $\epsilon > 0$, let $U_\epsilon$ be an open neighborhood of
$0$ on
which $|y(h_i)'| < \epsilon$ and $|y(h_i)(x)| < |x|\epsilon$.
Now, for two indices $i,j$ the composition has the form
\begin{align*}
h_i \circ h_j (x) & = x + y(h_j)(x) + y(h_i)(x + y(h_j)(x)) \\
& = x + y(h_j)(x) + y(h_i)(x) + O(\epsilon y(h_j)(x)).
\end{align*}
In particular, the composition deviates from $x + y(h_i)(x) + y(h_j)(x)$
by a
term which is small compared to $\max(y(h_i)(x),y(h_j)(x))$.

Now, choose some sequence of points $x_i \to 0$. For each $i$
define the map $v_i\co H \to \R^n$ where $v_i(h) = y(h)(x_i)$.
Let $w_i = \sup_{j \le m} |v_i(h_j)|$, and define $v_i'(h) = v_i/w_i$.
It follows that the functions $v_i'$ are uniformly bounded on each $h
\in H$,
and therefore there is some convergent subsequence. Moreover, by the
estimate above, on this subsequence, the maps $v_i'$ converge to a
homomorphism $v'\co H \to \R^n$. On the other hand, by construction,
there is some index $j$ such that $|v_i'(h_j)| = 1$. In particular, the
homomorphism $v'$ is nontrivial, and $H$ surjects onto a nontrivial
free abelian group, and we are done.
\end{proof}

\begin{lem}\label{locally_circular}
Let $G$ be a group of germs of orientation preserving
$C^1$ homeomorphisms of $\R^2$ at a fixed point $p$.
Then $G$ is circularly-orderable. Moreover, if the image $L(G)$ of $G$
in $\GL^+(2,\R)$ obtained by linearizing the action at $p$ is
left-orderable,
then $G$ itself is left-orderable.
\end{lem}
\begin{proof}
By the Thurston stability theorem, if $L\co G \to L(G)$ denotes the
homomorphism onto the linear part of $L(G)$, then the kernel $K(G)$
is locally indicable, and therefore by
Theorem~\ref{locally_surjective_LO},
$K(G)$ is LO. That is, we have a short exact sequence
$$0 \longrightarrow K(G) \longrightarrow G \longrightarrow L(G)
  \longrightarrow 0$$
where $K(G)$ is LO. If the image $L(G)$ is LO, then so is $G$ by
Lemma~\ref{ses_LO}.

Moreover, since $G$ is orientation preserving,
$L(G) < \GL^+(2,\R)$ where $\GL^+$ denotes the subgroup of $\GL$ with
positive determinant. There is a homomorphism from $\GL^+(2,\R)$ to
$\SL(2,\R)$ with kernel $\R^+$. We write this as a short exact sequence:
$$0 \longrightarrow \R^+ \longrightarrow \GL^+(2,\R) \longrightarrow
  \SL(1,\R) \longrightarrow 0.$$
The group $\SL(2,\R)$ double covers $\PSL(2,\R)$, and can be thought of as
the group of orientation-preserving transformations of the connected
double
cover of $\RP^1$ which are the pullback of projective transformations of
$\RP^1$ by $\PSL(2,\R)$. In particular, $\SL(2,\R)$ is a subgroup of
$\homeo^+(S^1)$,
and is therefore CO. It follows
by Lemma~\ref{ses_CO} that $\GL^+(2,\R)$ is CO, and therefore so is
$L(G)$.

By another application of Lemma~\ref{ses_CO}, $G$ is CO.
\end{proof}

\subsection{Groups which stabilize one or more compact regions}
The main point of this subsection is to prove
Theorem~\ref{fixes_components_circular},
which says that a group of $C^1$ orientation-preserving
homeomorphisms of the plane which fixes a compact set $K$ with connected
complement is circularly-orderable, and a group which fixes at least
two such sets
is left-orderable.

Now, suppose that a group $G$ stabilizes the disjoint compact connected
sets
$K_1,K_2,\dots$. For each $K_i$, exactly one complementary component is
unbounded. So without loss of
generality, we can fill in these bounded complementary
regions and assume $G$ stabilizes
disjoint compact connected sets $K_1, \dots$ with connected complement
for all $i$.

\begin{rmk}
{\rm A compact set $K_i$ which is an absolute neighborhood retract
satisfies the hypotheses of Alexander duality, and has
connected complement in $\R^2$ iff $H_1(K_i;\Z) = 0$. However, the
well-known example of the ``topologist's circle'' shows that in general,
the vanishing of homology is not sufficient to show that the complement
is connected.}
\end{rmk}

We show how each region $K_i$ gives rise to a CO group $G_i$ which is
naturally a quotient of $G$.

\begin{construct}\label{circular_components}
{\rm If $K_i$ consists of more that one point, the complement
of $K_i$ in the sphere $S^2 = \R^2 \cup \infty$ is conformally a disk, and
we let $G_i$ denote the image of $G$ in $\aut(\P(S^2 \backslash K_i))$.
Note that $G_i$ is CO, by Lemma~\ref{prime_ends_circular}.

If $K_i$ consists of a single point and $G$ acts $C^1$ near $K_i$, let
$G_i$ denote the germ of $G$ at $K_i$. By Lemma~\ref{locally_circular},
$G_i$ is CO, and is actually LO if the linear part $L(G)$ of $G$ at
$K_i$ is LO.}
\end{construct}

\begin{rmk}
{\rm Notice that there is a natural circular order on $G_i$ in the first
case, but {\em not} in the second. However, in the second case, the
{\em Euler class} of the circular order provided by
Lemma~\ref{locally_circular}
is just the circular Euler class of the linear part $L(G)$ of $G$
at $K_i$,
acting projectively on the unit tangent bundle, and is therefore natural.}
\end{rmk}

In particular, the set $K_i$ gives a homomorphism from $G$ to a product
$$G \longrightarrow \prod_i G_i$$
of CO groups. The product of any number of left-orderable groups
is left-orderable. This follows immediately from Lemma~\ref{ses_LO}
and Lemma~\ref{local_LO}. But
a product of CO groups is not necessarily CO. In this section we will show
that in the context above, if $G$ stabilizes two disjoint compact
connected sets
$K_1,K_2$ with connected complement, then the groups $G_1,G_2$ are
actually
both LO.

\begin{lem}\label{Euler_vanishes}
Let $G$ act by $C^1$ orientation-preserving homeomorphisms on the
plane. Suppose that
$G$ stabilizes two disjoint compact connected sets $K_1,K_2$ with
connected
complement. Then the groups $G_1,G_2$ provided by
Construction~\ref{circular_components}
are circularly-orderable, and the pullback of their Euler classes
to $G$ are zero.
\end{lem}
\begin{proof}
First suppose $K_1$ contains at least two points.

For a compact, connected set $K_1 \subset \R^2$ with connected complement,
the complement of $\R^2 \backslash K_1$ is
homeomorphic to an annulus $A_1$. Let $\til{A_1}$ denote the universal
cover of $A_1$, and denote by $\hat{G_1}$ the central extension of $G_1$
$$0 \longrightarrow \Z \longrightarrow \hat{G_1} \longrightarrow G_1
  \longrightarrow 0.$$
This requires some explanation. The Riemann surface $\til{A_1}$
is noncompact,
and topologically is an open disk. Under a uniformizing map, there is
a unique
point $p \in S^1$ such that every proper ray $r$ in $\til{A_1}$ which
projects to an unbounded proper ray in $A_1$ is asymptotic to $p$. If
we cut
$S^1$ at $p$, we get a copy of $\R$, which by abuse of notation we refer
to as
the set of prime ends $\P(\til{A_1})$

Since by hypothesis $K_1$ consists of more than one point,
there is a natural map $\P(\til{A_1}) \to \P(S^2\backslash K_1)$ which
is just a
covering map $\R \to S^1$ under the identification of $\P(S^2 \backslash
K_1)$
with $\P(D) = S^1$ by a uniformizing map.
Then the group $\hat{G_1}$ is just the usual preimage
of a subgroup of $\homeo^+(S^1)$ in ${\thomeo^+(S^1)}$.

In particular, the class of this extension is the Euler class of the
circular
ordering on $G_i$. In order to show that the pullback of this Euler
class to
$G$ is trivial, it suffices to show that the restriction
homomorphism $\text{res}\co G \to \homeo^+(A_1)$ lifts to
the covering space $\hat{\text{res}}\co  G \to \homeo^+(\til{A_1})$.

Let $\hat{K_2}$ be a lift of $K_2$ to $\til{A_1}$. Then for each element
$g \in G$, there is a unique lift $\hat{\text{res}}\co g \to
\homeo^+(\til{A_1})$ which
stabilizes $\hat{K_2}$. By uniqueness, this defines the desired section.

The case that $K_1$ consists of a single point is very similar. The
complement
$\R^2 \backslash K_1$ is again an annulus, and there is
A covering $\til{A_1} \to A_1$. There is a central extension
$$0 \longrightarrow \Z \longrightarrow \hat{G_1} \longrightarrow G_1
  \longrightarrow 0$$
where $G_1$ is the germ of $G$ at $K_1$, and $\hat{G_1}$ is the
corresponding germ of the
preimage of $G$ in $\homeo^+(\til{A_1})$. Again, choosing a lift
$\hat{K_2}$ of $K_2$ to $\til{A_1}$ determines a unique section
$G \to \homeo^+(\til{A_1})$ whose germ is $\hat{G_1}$. So the pullback of
the Euler class to $G$ is trivial in this case too.
\end{proof}

We now establish a technical lemma about groups which act in a
$C^1$ fashion and have indiscrete fixed point set.

\begin{lem}\label{not_isolated_fixed_points}
Suppose $G$ acts faithfully by $C^1$ orientation-preserving
homeomorphisms of the
plane, and suppose that the fixed point set $\fix(G)$ is not discrete.
Then $G$ is left-orderable.
\end{lem}
\begin{proof}
Let $p \in \fix(G)$ in the frontier of $\fix(G)$
be a limit point of distinct $p_i \in \fix(G)$. Then
the $p_i$ contain a subsequence which approach $p$ radially along
some vector
$v$. It follows that the linearization of $G$ at $p$, denoted $L(G)$,
fixes the
vector $v$, and therefore the image is contained in some conjugate of the
affine group of the line $\aff^+(\R)$. The affine group of the line is an
extension of $\R$ by $\R$:
$$0 \longrightarrow \R \longrightarrow \aff^+(\R) \longrightarrow \R
  \longrightarrow 0$$
where the homomorphism to $\R$ is given by the log of the dilation,
and the kernel is the subgroup of
translations. It follows that $\aff^+(\R)$
is left-orderable by Lemma~\ref{ses_LO}.

It follows from Lemma~\ref{locally_circular} that the germ of $G$ at $p$
is LO.
In particular, if $H$ is any finitely-generated subgroup of $G$, then
we can
apply this reasoning to a limit point in the frontier of $\fix(H)$ and see
that the germ of $H$ at this limit point is left-orderable. In particular,
$H$ is
LO--surjective. Since $H$ was arbitrary, by
Theorem~\ref{locally_surjective_LO}
the group $G$ is LO as required.
\end{proof}

We now have all the tools to prove the main result of this subsection:

\begin{thm}\label{fixes_components_circular}
Let $G$ act faithfully by $C^1$ orientation-preserving homeomorphisms
of the plane.
Suppose that $G$ stabilizes disjoint compact connected sets $\lbrace
K_i \rbrace$ such
that each $K_i$ has connected complement. Then $G$ is CO. If there is
more than
one region $K_i$, then $G$ is LO.
\end{thm}
\begin{proof}
By Lemma~\ref{Euler_vanishes} there is a homomorphism
$$G \longrightarrow G_i$$
for each $i$, where the image group is CO, and the pullback of the Euler
class of
this circular ordering to $G$ is trivial if there is more than one
region $K_i$.

Let $K$ denote the kernel of this homomorphism. For each region $K_i$ with
more than one point, Lemma~\ref{action_faithful} implies that $K$ fixes
$\partial K_i$ pointwise. For each region $K_i$ consisting of a single
point,
every element of $K$ fixes some neighborhood of $K_i$ pointwise. It
follows
that for each finitely-generated subgroup $H$ of $K$ the set of fixed
points
of $H$ is not isolated. In particular, by
Lemma~\ref{not_isolated_fixed_points},
every finitely-generated subgroup of $K$ is LO, and therefore $K$ is LO by
Lemma~\ref{local_LO}.

By Lemma~\ref{ses_CO}, $G$ is CO. Moreover, the Euler class of this
circular
ordering is the pullback of the Euler class of $G_i$ under the
homomorphism $G \to G_i$. In particular, if there is more than one
region $K_i$,
the Euler class of the circular ordering on $G$ is zero, and therefore $G$
is LO, as claimed.
\end{proof}

\subsection{Totally disconnected planar sets}

If $G$ acts faithfully on $\R^2$ and merely {\em permutes} a collection of
disjoint compact regions $\lbrace K_i \rbrace$ we can study the image
of $G$ in the
group of homotopy classes of {\em all} orientation preserving
homeomorphisms
of $\R^2$ which permute the $\lbrace K_i \rbrace$. We will call such
a group
a generalized braid group. The following construction is a convenient
simplification.

\begin{construct}\label{fill_in_and_shrink}
{\rm Suppose $G$ acts faithfully on $\R^2$ and preserves some bounded set $K$.
Take the closure of $K$, fill in bounded complementary regions, and label
the components of the result as $\lbrace K_i \rbrace$.

Since the $K_i$ are obtained naturally from $K$,
$G$ permutes the components $\lbrace K_i \rbrace$. Observe that the
union of
the $K_i$ is compact, and that the $K_i$ are disjoint, compact, and
have connected complement.
Then we can define an equivalence relation on
$\R^2$ by quotienting each $K_i$ to a point.}
\end{construct}

\begin{defn}
{\rm A {\em decomposition} $\mathscr{G}$ of a space is a partition into
compact subsets. A
decomposition is {\em upper semi-continuous} if for every decomposition
element
$\zeta \in \mathscr{G}$ and every open set $U$ with $\zeta \subset U$,
there exists
an open set $V \subset U$ with $\zeta \subset V$ such that every $\zeta'
\in \mathscr{G}$
with $\zeta' \cap V \ne \emptyset$ has $\zeta' \subset U$. The
decomposition is
{\em monotone} if its elements are connected.}
\end{defn}

It is clear that the equivalence relation in
Construction~\ref{fill_in_and_shrink}
is monotone and upper semi-continuous.

The following theorem is proved by R. L. Moore in \cite{Moore}:
\begin{thm}[Moore]\label{decomposition_theorem}
Let $\mathscr{G}$ be a monotone upper semi-continuous decomposition of
a topological
sphere $S$ such that no decomposition element separates $S$. Then the
quotient space
obtained by collapsing each decomposition element to a point is a
topological sphere.
\end{thm}

\begin{defn}
{\rm Suppose $G$ acts on a space $X$ and preserves a monotone equivalence
relation
$\sim$. Then the action of $G$ descends to an action on $X/\sim$. We say
the two actions are {\em monotone equivalent}.}
\end{defn}

\begin{lem}\label{totally_disconnected_quotient}
Let $G$ act on $\R^2$ in such a way as to permute a union of disjoint
compact
regions. Then this action is monotone equivalent to an action which leaves
invariant a compact, totally disconnected, $G$--invariant subset.
\end{lem}
\begin{proof}
The fact that $\R^2/\sim$ is homeomorphic to $\R^2$ follows from Moore's
Theorem~\ref{decomposition_theorem}. Since
Construction~\ref{fill_in_and_shrink}
is natural, the equivalence relation it defines is $G$--invariant, and
therefore the action of $G$ descends to the quotient. It is clear from the
construction that $C$ is totally disconnected.
\end{proof}

Since $G$ acts on $\R^2$ and permutes $C$, we can study a {\em minimal
subset} of $C$; that is, a closed invariant nonempty subset of $\R^2$
which is
minimal with respect to these properties. Amongst all invariant subsets,
a minimal subset $C$ can be characterized by
the property that every $G$ orbit contained in $C$ is dense in $C$.
Such a subset will be either
finite or {\em perfect}. Here, a set is {\em perfect} if every point is
a limit
point.

The Moore--Kline theorem gives a useful ``normal form'' for totally
disconnected compact planar sets. This theorem is proved in
\cite{Moore_Kline}:

\begin{thm}[Moore, Kline]\label{disconnected_in_arc}
Let $C$ be a totally disconnected compact set in $\R^2$. Then there is a
homeomorphism of $\R^2$ to itself which takes $C$ to a closed subset
of the
arc $[0,1]$ contained in the $x$--axis.
\end{thm}

We will use this theorem in the next subsection. Note that this theorem
implies
that a perfect, totally disconnected subset of $\R^2$ is a tame Cantor
set.

\subsection{Generalized braid groups}

\begin{defn}
{\rm For a surface $S$ (possibly with boundary $\partial S$) the {\em relative
orientation-preserving mapping class group},
denoted $\MCG^+(S,\partial S)$, is the group
of homotopy classes of orientation-preserving homeomorphisms of $S$ to
itself which are fixed on $\partial S$. We generally just call this the
mapping class group for short.

Let $C$ be a compact, totally disconnected subset of the open unit disk
$D$. Then the {\em generalized braid group} of $C$, denoted $B_C$,
is the relative mapping
class group $\MCG^+(\overline{D}\backslash C,\partial \overline{D})$ and
the {\em planar generalized braid group} of $C$, denoted $B_C'$, is the
mapping class group
$\MCG^+(\R^2\backslash C)$.}
\end{defn}

By extending a homeomorphism of $D$ fixed on
$\partial \overline{D}$ by the identity, we see that there is a natural
homomorphism $B_C \to B_C'$. This gives rise to a central extension
$$0 \longrightarrow \Z \longrightarrow B_C \longrightarrow B_C'
  \longrightarrow 0$$
where the generator of $\Z$ is the (positive) Dehn twist of $D^2$ along
a boundary
parallel loop.

If $C$ is a finite set with $n$ points, then $B_n$ is the usual braid
group,
and $B_n'$ is the usual quotient of $B_n$ by its center. A reference for
the theory of braid groups, with particular reference to questions
of orderability
or circular orderability, is \cite{Dehorn_book}.

We come to the first main result of this paper.

The proof of this theorem follows the same general approach as Thurston's
proof of
the orderability of the usual (finitely-generated) braid groups
\cite{Short_Weiss},
but several technical complications arise because the surfaces $D^2
\backslash C$ and
$\R^2 \backslash C$ are not of finite type.

\begin{mapping_class_thm}\label{mapping_class_thm}
Let $C$ be a compact, totally disconnected subset of the open unit
disk $D$.
Then $B_C'$ is circularly-orderable, and $B_C$ is left-orderable.
\end{mapping_class_thm}
\begin{proof}
If $C$ consists of a single point, both groups are trivial. If $C$
consists of two points, $B_C' = \Z/2\Z$ and $B_C = \Z$, so the theorem
is satisfied in that case. So without loss of generality we
assume that $C$ contains at least $3$ points.

First we show that $B_C'$ is circularly-orderable.
Let $S = S^2 \backslash C$.
Then $S$ is a hyperbolic surface, and we
can identify the universal cover $\til{S} = \H^2$. We work in the unit
disk model of $\H^2$, so that $\til{S}$ can be compactified by adding
the circle at infinity $S^1_\infty$. There is a short exact sequence
$$0 \longrightarrow \pi_1(S) \longrightarrow \hat{\MCG(S)} \longrightarrow
  \MCG(S) \longrightarrow 0.$$
Here $\hat{\MCG(S)}$ denotes the group of homotopy classes of
homeomorphisms
of $\til{S}$ which commute with the action of $\pi_1(S)$, and cover
homeomorphisms
of $S$.

It turns out that there is a naturally defined injective
homomorphism $\sigma$ from the mapping class group $\MCG(\R^2\backslash
C)$
to $\hat{\MCG(S)}$. We define $\sigma$ as follows.
Let $[\phi] \in \MCG(\R^2\backslash C)$ denote a typical
element, and let $\phi$ be a representative. Then $\phi$ extends
continuously
to a homeomorphism $\phi_S$ of $S$ which fixes the point $\infty$. Choose
a lift
$\hat{\infty}$ of $\infty$ in $\til{S}$. Then there is a unique lift of
$\phi_S$ to $\til{S}$ which fixes $\hat{\infty}$; call this lift
$\phi_{\til{S}}$.
Then define $\sigma$ by
$$\sigma([\phi]) = [\phi_{\til{S}}].$$
Since the construction of $\phi_{\til{S}}$ depends continuously on
$\phi$, it
follows that the class $[\phi_{\til{S}}]$ is well-defined.

We show that it is injective. If
$[\phi_{\til{S}}] = [\psi_{\til{S}}]$ for some $\phi,\psi$ then there
is a homotopy $\Psi_{\til{S}}\co \til{S} \times I \to \til{S}$ with
$\Psi_{\til{S}}(\cdot,0) = \psi_{\til{S}}$ and
$\Psi_{\til{S}}(\cdot,1) = \phi_{\til{S}}$ which at every point is
$\pi_1(S)$ equivariant. It follows that $\Psi_{\til{S}}$ descends to a
map $\Psi_S\co S \times I \to S$ which is a homotopy between $\psi_S$ and
$\phi_S$. Moreover, under the track of this homotopy, the point $\infty$
moves in a homotopically inessential loop, since by construction,
this loop lifts to $\til{S}$. So we can find another homotopy $\Psi_S'$,
homotopic
to $\Psi_S$, which fixes $\infty$. It follows that $\Psi_S'$ restricts to
$\Psi\co \R^2\backslash C \times I \to \R^2 \backslash C$ which is
a homotopy from $\psi$ to $\phi$. This proves that the map is an
injection.

In fact, it is not hard to show that the map $\sigma$ is actually an
isomorphism, but this is unnecessary for our purposes.

We now show that the group $\hat{\MCG(S)}$ is circular.
First of all we construct an action of $\hat{\MCG(S)}$ on a
circularly-ordered set.

The conformal structure that $S$ inherits as an open submanifold of
$\R^2$ determines a natural complete hyperbolic structure on $S$.
This gives an identification of the universal cover $\til{S}$ with
the hyperbolic plane $\H^2$, which can be compactified by its ideal
boundary which is a circle $S^1_\infty$. An orientation on $S$ determines
one on $\til{S}$ and therefore also on $S^1_\infty$.
Since the hyperbolic structure on $S$ is complete,
the limit set $\Lambda$ of $\pi_1(S)$ is the entire
circle $S^1_\infty$. A dense subset of $S^1_\infty$ consists of
endpoints of
axes of hyperbolic translations $\alpha \in \pi_1(S)$. Notice that since
$\pi_1(S)$ is discrete, the stabilizer of a point in $\pi_1(S)$ is cyclic,
so distinct axes have distinct endpoints. Such endpoints are in bijective
correspondence with the set $\E$ whose elements are
maximal hyperbolic cyclic subgroups
$\langle \alpha \rangle$ of $\pi_1(S)$ together with a choice of generator
for $H_1(\langle \alpha \rangle;\Z)$. It follows that $\E$ inherits
a circular
order by its bijection with a subset of the circle $S^1_\infty$.

The group $\hat{\MCG(S)}$ induces
an automorphism of $\pi_1(S)$ by conjugation, and therefore induces an
action on the set $\E$. We verify that this action preserves the
cyclic order.

To see this, let $\alpha,\beta$ be hyperbolic elements of $\pi_1(S)$ with
axes $l_\alpha, l_\beta$. Let $[\phi] \in \hat{\MCG(S)}$, and let
$\phi$ be a representative. The images of the axes
$\phi(l_\alpha),\phi(l_\beta)$
are periodic under the elements
$[\phi]_*(\alpha) = [\phi]\alpha[\phi]^{-1}$ and
$[\phi]_*(\beta) = [\phi]\beta[\phi]^{-1}$ in $\pi_1(S)$
respectively, so they converge to well-defined end points in
$S^1_\infty$, which
are exactly the endpoints of the axes of
$l_{[\phi]_*(\alpha)},l_{[\phi]_*(\beta)}$.
Moreover, they are either disjoint or cross transversely in exactly one
point, according to whether $l_\alpha,l_\beta$ are disjoint or cross
transversely
in exactly one point respectively. It follows that $l_\alpha,l_\beta$
cross
iff $l_{[\phi]_*(\alpha)}, l_{[\phi]_*(\beta)}$ cross. See
Figure~\ref{axes}.

\begin{figure}[ht]
\centerline{\relabelbox\small \epsfxsize 4.0truein
\epsfbox{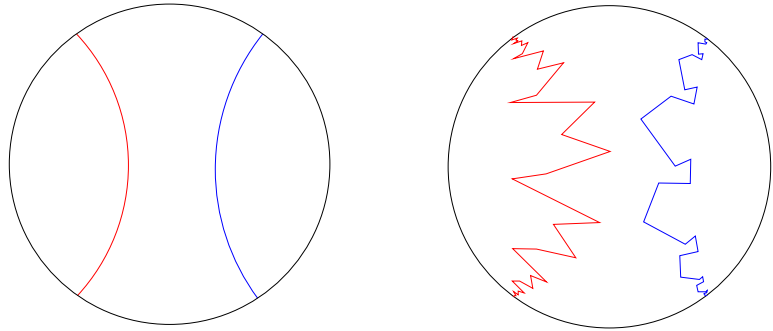}
\endrelabelbox}
\caption{The axes $l_\alpha,l_\beta$ are disjoint iff their images
under $\phi$ are disjoint. But this is true iff the axes
$l_{[\phi]_*(\alpha)},l_{[\phi]_*(\beta)}$ are disjoint. This defines a
natural action of $\hat{\MCG(S)}$ on the circularly-ordered set of
endpoints of axes.}
\label{axes}
\end{figure}

The circular ordering on the endpoints
of $l_\alpha,l_\beta$ is determined by this crossing information, the
orientation of the axes, and if the axes do not intersect, the
relationship
between the orientation of the axes and the orientation of the
subsurface of
$\til{S}$ that they both cobound. It follows that the circular order of
the endpoints
of the $l_\alpha,l_\beta$ is taken to the circular order of the endpoints
of the axes
$l_{[\phi]_*(\alpha)},l_{[\phi]_*(\beta)}$ under $[\phi]$, and therefore
the action of $\hat{\MCG(S)}$ preserves the circular order on $\E$.

It remains to show this action is faithful.
One way to do this is to find
a maximal collection of pairwise disjoint geodesic loops
$\Gamma = \lbrace \gamma_i \rbrace$
in $S$. We want to be slightly careful about the choice of $\Gamma$; we
want it to be maximal, and also {\em closed} as a subset of $S$. In
particular,
the geodesics $\gamma_i$ should not accumulate. We show how to choose
such a
collection.

By the Moore--Kline Theorem~\ref{disconnected_in_arc} we can assume, after
a homeomorphism of $\R^2$, that the set $C$ is contained as a closed,
totally disconnected subset of a horizontal arc $I$ in the $x$--axis. In
particular, the
complement $I \backslash C$ consists of a countable collection of open
intervals $I_i$.
For each connected complementary open interval $I_i$ we let $p_i$
denote the
midpoint. We inductively produce a maximal set of round circles $S_k$
which are
symmetric in the $x$--axis, and intersect the $x$ axis in points $p_i,p_j$
subject
to the constraint that the circles are
non-intersecting, and the $k$th circle is chosen to be the biggest
circle with
diameter $\le 2^{-k}$. Note that a maximal collection of such
circles has the property that each complementary domain contains at most
finitely many unmatched $p_i$.

Our initial guess for $\Gamma$ is the set of geodesics isotopic to $S_i$.
Note that some $S_i$ might be boundary parallel; we throw these out and
just look at the ones with geodesic representatives.
This is not quite maximal, but each complementary region is of finite
type,
and there are only finitely many complementary regions whose diameter
(in the Euclidean metric in $\R^2$) is $> \epsilon$ for any $\epsilon$.
So we just take a maximal (finite) collection of geodesic loops in
each such
complementary region, and add these loops to $\Gamma$. Notice that this
union of geodesic loops is {\em closed}, doesn't accumulate anywhere
(in $S$)
and is maximal.

Now, by hypothesis, $[\phi]$ preserves every hyperbolic cyclic subgroup
$\langle \alpha \rangle$ of $\pi_1(S)$; in particular, the projection
$\phi_S$ of a representative $\phi$ to $S$ preserves every oriented
free homotopy
class of loop which is not boundary parallel. So the collection $\Gamma$
is
taken to a collection of curves which are curvewise isotopic to
$\Gamma$. By
straightening these geodesics inductively, we can homotope $\phi_S$ to
$\phi_S'$ which fixes $\Gamma$ pointwise. Since $\Gamma$ was maximal,
every
complementary surface to $\Gamma$ is either a thrice punctured sphere,
a twice punctured disk, a once punctured annulus, or a pair of pants.
An automorphism of one of these subsurfaces which preserves every
hyperbolic
cyclic subgroup is isotopic to the identity, so $\phi_S$ is actually
isotopic
to the identity, and therefore $[\phi]$ is in $\pi_1(S)$. Now, an element
of $\pi_1(S)$ which, under conjugation,
takes a hyperbolic cyclic subgroup to itself in an
orientation preserving way must actually commute with that subgroup. But
this contradicts the fact that $\pi_1(S)$ is nonelementary.

The conclusion is that $\hat{\MCG(S)}$ acts faithfully in an order
preserving
way on a circularly-ordered set, and therefore $\hat{\MCG(S)}$ is
circularly-ordered, by Theorem~\ref{action_is_order}. Since $\sigma\co
B_C' \to
\hat{\MCG(S)}$ is injective, $B_C'$ is circularly-ordered too.

Proving that $B_C$ is left-orderable is similar; let $p \in \R^2\backslash
D$ be some point.
Define $S$ to be obtained from $\R^2 \backslash (C \cup p)$ by adding a
point to compactify $\R^2$, and then removing $(C \cup p)$. Exactly as in
the previous case, there is an injective homomorphism
from $B_C$ to $\hat{\MCG(S)}$, and an invariant circular order on
$\hat{\MCG(S)}$.
Moreover, the image of $B_C$ stabilizes the parabolic subgroup
corresponding
to a small loop around $p$, since representatives of $B_C$ act trivially
there.
In particular, the image of $B_C$ stabilizes a point in $S^1_\infty$, and
therefore it is left-orderable, as required.
\end{proof}

\subsection{Groups with bounded orbits}

We are now in a position to prove the following:

\begin{bounded_orbit_thm}
Let $G$ be a group of orientation preserving $C^1$ homeomorphisms of
$\R^2$ with a bounded orbit. Then $G$ is circularly-orderable.
\end{bounded_orbit_thm}
\begin{proof}
Let $p$ be a point with a bounded orbit, and consider the
closure of the orbit of $p$. This is a compact union
of disjoint compact components which is permuted by $G$.

We apply
Construction~\ref{fill_in_and_shrink} to fill in the sets with
disconnected complement
to get a new collection of connected sets $\lbrace K_i \rbrace$ which are
compact with connected complements, whose union is compact,
and whose components are permuted by $G$.
We then apply Lemma~\ref{totally_disconnected_quotient} to get
a monotone equivalent action of $G$ on $\R^2$ which leaves invariant
a totally disconnected compact set $C$.

By theorem A, the image of $G$ in the mapping class group
of $\R^2\backslash C$
is circularly-orderable. Denote the kernel by $K$. Then $K$
fixes $C$ pointwise, and therefore with respect to the original action,
$K$ preserves the components $\lbrace K_i \rbrace$ componentwise.

By Theorem~\ref{fixes_components_circular}, $K$ is circularly-orderable,
and is left-orderable if there are at least two components of $K_i$.
If there is only one component $K_i$, the quotient $C$ is a single point,
and the mapping class group of $\R^2\backslash C$
is trivial, in which case $G=K$.
So either $G$ is circularly-orderable, or else it maps to a
circularly-orderable group with left-orderable
kernel. By Lemma~\ref{ses_CO}, $G$ is circularly-orderable, as claimed.
\end{proof}

\section{The Euler class and smoothness}

The purpose of this section is to give a complete classification of which
Euler classes can arise for orientation-preserving
actions of oriented surface groups on $\R^2$ in every degree of
smoothness.

One main result is that the Euler class of the group of $C^\infty$
orientation-preser\-ving diffeomorphisms of the plane is an unbounded class,
in stark contrast with the case of $\homeo^+(S^1)$ and the Milnor--Wood
inequality.  In particular, if $S$ is a closed surface of genus at
least $2$, there is a $C^\infty$ action of $\pi_1(S)$ on the plane
with arbitrary Euler class.  The other main result is that the genus
one case is more rigid: we show that an orientation-preserving $C^1$
action of $\Z \oplus \Z$ on the plane has vanishing Euler class.

Using these facts, we prove the existence of a finitely-generated
torsion-free group
of homeomorphisms of the plane which is not circularly-orderable.

\subsection{The Euler class and winding number}

We give two descriptions of the Euler class of a group $G$ acting
on $\R^2$
by a representation $\rho\co G \to \homeo^+(\R^2)$. The first description
is
more algebraic.

\begin{construct}\label{algebraic_class}
{\rm Let $G_\infty$ denote the {\em germ of $G$ at $\infty$}.
Let $A$ be a punctured disk neighborhood of $\infty$, and let $\til{A}$
be the universal cover of $A$. Then there is a central extension
$$0 \longrightarrow \Z \longrightarrow \hat{G}_\infty \longrightarrow
  G_\infty \longrightarrow 0$$
where $\hat{G}_\infty$ denotes the subgroup of periodic germs of
homeomorphisms
of $\til{A}$. The class of this extension pulls back by the natural
restriction homomorphism $G \to G_\infty$ to give the
Euler class $\rho^*[e] \in H^2(G;\Z)$.}
\end{construct}

To see that this is the Euler class, observe that it defines a non-split
extension, and therefore determines a non-trivial element of
$H^2(\homeo^+(\R^2);\Z)$.
By evaluating this class on some suitable homology cycles (e.g.
Example~\ref{Bestvina_example} in the sequel) one sees
that it is a primitive element, and therefore it is equal
to the generator of $H^2(\homeo^+(\R^2);\Z)$, as claimed.

The next description is geometric, in terms of a winding number. Notice
that for this
definition to make sense, we must assume the action is at least $C^1$.

\begin{construct}\label{geometric_class}
{\rm Now we give a more geometric construction of the Euler class, in terms of
how it pairs with elements of $H_2(G;\Z)$. Assume the action is at
least $C^1$.
It is a fact that for any group $G$, an integral class $\sigma \in
H_2(G;\Z)$
can be realized as the image of a map from a closed orientable surface
$S$ of
genus $\ge 1$ into a $K(G,1)$
$$f\co S \longrightarrow K(G,1).$$
Given a representation $\rho$, we get an induced action
$$\rho\co \pi_1(S) \longrightarrow \homeo^+(\R^2).$$
Let $g$ be the genus of $S$, and suppose that we obtain
$S$ by gluing the sides of a $4g$--sided polygon $P_g$ in pairs in the
usual way. Let $d_\rho\co \til{S} \to \R^2$ be an equivariant $C^1$ map,
called a {\em developing map} and
look at the restriction $$d_\rho\co  \partial P_g \longrightarrow \R^2$$
where $P_g$ now denotes a fundamental domain, by abuse of notation.
If necessary, we can assume that $d_\rho$ is an immersion on edges of
$\partial P_g$.

Let $\delta\co  S^1 \to \R^2$ be obtained from $d_\rho$ by smoothing the
corners at each vertex of $\partial P_g$. Then the
Euler class of $\rho$, evaluated on $[S]$, can be calculated by the
formula
$$\rho^*([e])([S]) = \index(\delta) + 1 - 2g$$
where $\index(\delta)$ denotes the usual winding number of the tangent
vector $T\delta$ around $S^1$. See Figure~\ref{euler_index} for
an example of a hypothetical action of $\Z \oplus \Z$ with
Euler number $1$. In the figure, $\alpha,\beta$ denote generators
of $\Z \oplus \Z$.

\begin{figure}[ht]
\centerline{\relabelbox\small \epsfxsize 2.0truein
\epsfbox{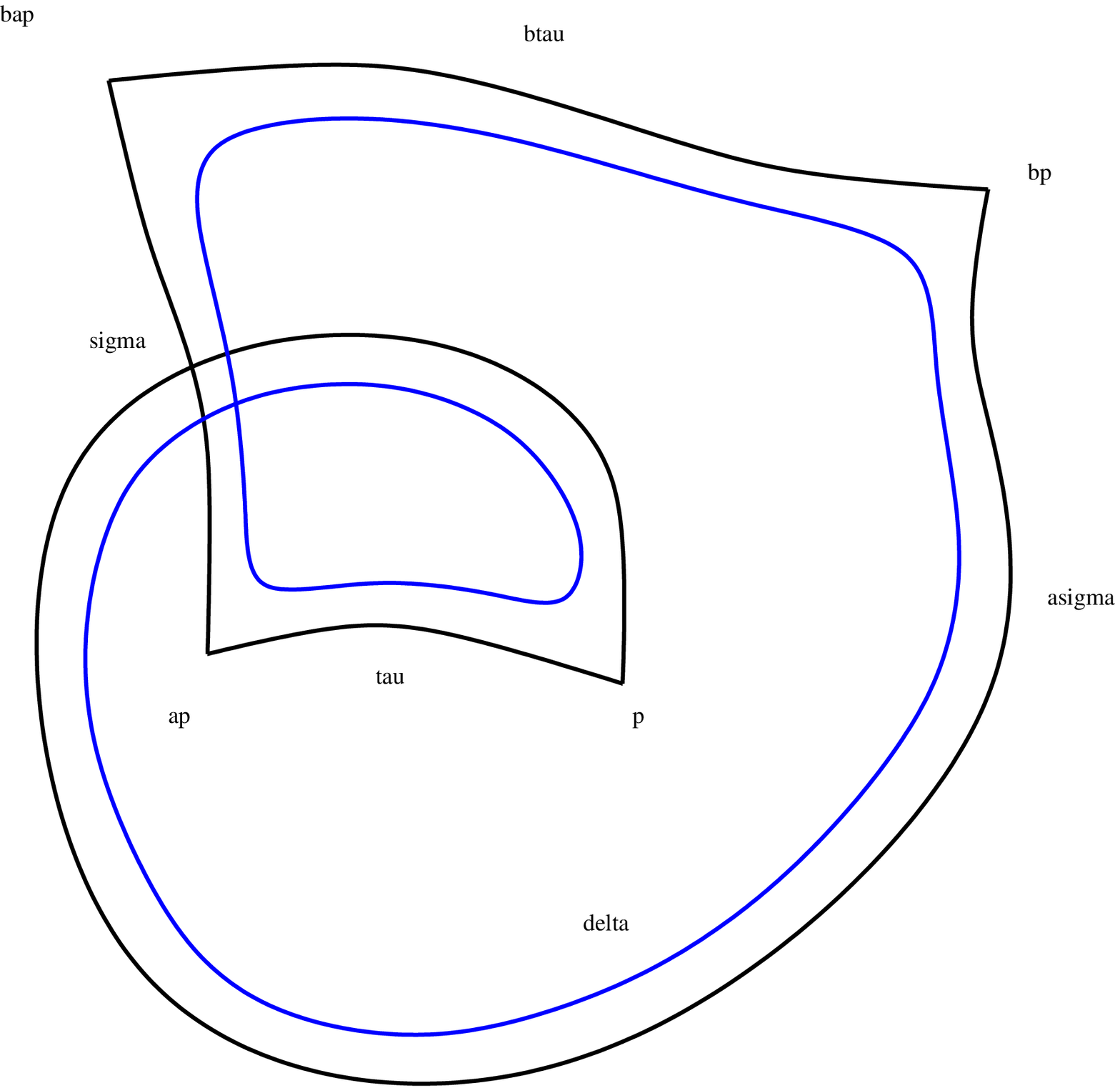}
\relabel{p}{$p$}
\relabel{ap}{$\alpha(p)$}
\relabel{bp}{$\beta(p)$}
\relabel{bap}{$\beta\alpha(p)$}
\relabel{delta}{$\delta$}
\relabel{tau}{$\tau$}
\relabel{sigma}{$\sigma$}
\relabel{btau}{$\beta(\tau)$}
\relabel{asigma}{$\alpha^{-1}(\sigma)$}
\endrelabelbox}
\caption{This hypothetical action of $\Z \oplus \Z$ satisfies
$\text{index}(\delta)=2$, so $\rho^*([e])([\Sigma])=1$.}
\label{euler_index}
\end{figure}

In practice we would like to work with developing maps which are as
degenerate as
possible, in order to reduce the combinatorial complexity of the
calculation.
So to calculate the ``correct'' smoothing along a degenerate edge, we must
perturb the picture near a degenerate edge in an equivariant way;
i.e. the perturbations
must be combinatorially equivalent along degenerate edges which are
paired by
the edge pairing translations.

The justification for this formula is as follows.
Associated to a representation
$\rho\co G \to \homeo^+(\R^2)$ there is an associated {\em foliated $\R^2$
bundle} $E_\rho$ over $BG$. Here $BG$ is the classifying space for $G$
as a discrete group; that is, it is a $K(G,1)$.
It is defined as the quotient of the product
$\R^2 \times EG$ by the action of $G$, where
$EG$ is the universal cover of $BG$, and
$$\gamma(p,q) = (\rho(\gamma)(p),\gamma(q))$$
By basic obstruction theory, the Euler class $\rho^*([e])$
is the obstruction to trivializing
$E_\rho$ as an $\R^2$ bundle over the $2$--skeleton of $BG$.
See for example \cite{Husemoller} for details.

If $\Sigma$ is a closed, orientable surface of genus $g$, and
$G = \pi_1(\Sigma)$, then $BG = \Sigma$. Since $\R^2$ is
contractible, we can always find a {\em section} $\sigma\co \Sigma \to
E_\rho$
of $E_\rho$ over $\Sigma$.

Suppose $\rho$ is $C^1$. Then $E_\rho$ is naturally a $C^1$ manifold,
and we can
choose $\sigma$ to be smooth. The pullback of the unit tangent bundle
in the fiber
direction by $\sigma$ defines an orientable circle bundle $E'_\rho$
over $\Sigma$.
A trivialization of $E'_\rho$ defines a trivialization of $E_\rho$
by exponentiation,
and conversely, a (smooth) trivialization of $E_\rho$ defines a
trivialization of
$E_\rho$ by restriction. Since $\diffeo^+(\R^2)$ and $\homeo^+(\R^2)$
are homotopic
as topological groups, $E_\rho$ can be trivialized iff it can be smoothly
trivialized.

It follows that the Euler class of $\rho$ is the Euler class of the
circle bundle
$E'_\rho$. The section $\sigma$ defines a developing map $d_\rho$. If
necessary,
perturb $d_\rho$ to be an immersion on each edge of $P_g$. This immersion
gives
a trivialization of $E'_\rho$ along the pushoff (into $P_g$) of each
edge of $\partial P_g$. This trivialization on either side of an edge of
the $1$--skeleton disagrees along the edge, and therefore the
trivialization
can be extended across each edge of the $1$--skeleton of $S$ with a single
saddle singularity, and across the vertex with a single sink singularity.
This section extends non-singularly over a fundamental domain $P_g$
iff the index of the smoothed curve $d_\rho(\delta)$ is $0$.
It follows that the Euler class, evaluated on $\Sigma$, is
$\index(\delta) + 1 - 2g$. The formula follows.}
\end{construct}

\begin{figure}[ht]
\centerline{\relabelbox\small \epsfxsize 3.0truein
\epsfbox{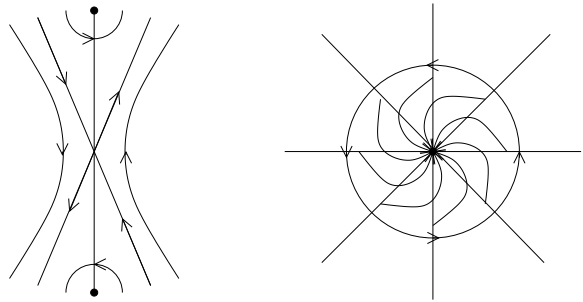}
\endrelabelbox}
\caption{The trivialization can be extended over each edge with a saddle
singularity, and over the vertex with a sink singularity.}
\label{smooth_corner}
\end{figure}

See Figure~\ref{smooth_corner} for a picture of the trivialization in
a neighborhood of
a typical edge and vertex in the genus $2$ case.

\subsection{Unboundedness of the Euler class}

In this subsection we construct $C^\infty$ actions on the plane
of $\pi_1(S)$, where $S$ is the closed surface of genus $2$, with
arbitrary
Euler class.

We begin with the following example, due to Bestvina:

\begin{exa}[Bestvina]\label{Bestvina_example}
{\rm Let $C_i$ be the round circle in $\R^2$ centered at the origin,
with radius $2^i$, for $i \in \Z$. That is, the $C_i$ are a bi-infinite
nested sequence of circles about $0$.
Let $\alpha$ be the composition of a positive Dehn twist in each circle
$C_i$.

More precisely, let $A_i$ be the annulus whose boundary
components are the circle of radius $0.99 \times 2^i$
and the circle of radius $1.01 \times 2^i$. Let $\psi\co [0,1] \to [0,1]$
be a
smooth homeomorphism, infinitely tangent to a constant map at $0$ and $1$.
Then let $\alpha$ be the identity outside the disjoint union of the
$A_i$, and
let its restriction to the circle of radius $(0.99 + t \times 0.02)
\times 2^i$
be a rotation through angle $2\pi \psi(t)$ for $t \in [0,1]$. In
particular, note that $\alpha$
is smooth away from the origin.

Let $\beta$ be the dilation centered at $0$ sending $z \to 2z$. Then
the commutator $[\alpha,\beta] = \id$, and it is easy to see that
$$\langle \alpha,\beta \rangle \cong \Z \oplus \Z.$$
On the other hand, the image of the boundary of the polygon $\partial P_1$
under $d_\rho$ from Construction~\ref{geometric_class} is isotopic to
the configuration in Figure~\ref{euler_index}. In particular, the
Euler number
of this $\Z \oplus \Z$ action is $1$.

If $G_i$ is a subgroup of $\Z \oplus \Z$ of index $i$, then $G_i$ is also
isomorphic to $\Z \oplus \Z$, and the induced action has Euler number $i$.
In particular, the Euler class $[e] \in H^2(\homeo^+(\R^2);\Z)$ is
{\em unbounded}.}
\end{exa}

The action of $\Z \oplus \Z$ is not $C^1$ at the origin. Bestvina
posed the
question of whether the Euler class was unbounded for $C^1$ actions. We
show
that the answer to this question is positive, and in fact show that
the Euler class is unbounded for $C^\infty$ actions.

\begin{Euler_class_thm}
For each integer $i$, there is a $C^\infty$
action $$\rho_i\co \pi_1(S) \to \diffeo^+(\R^2)$$
where $S$ denotes the closed surface of genus $2$, satisfying
$$\rho_i^*([e])([S]) = i.$$
In particular, the Euler class $[e] \in H^2(\diffeo^+(\R^2);\Z)$ is
{\em unbounded}.
\end{Euler_class_thm}
\begin{proof}
As in Example~\ref{Bestvina_example}, let $C_j$ be the round circle in
$\R^2$ with radius $2^j$ centered at the origin.
Let $\alpha$ be the composition of a positive
Dehn twist in each circle $C_j$ for $j \ge 0$, and let $\beta$ be
the dilation centered at $0$ sending $z \to 2z$. Then $\alpha,\beta$
are both
$C^\infty$, and the commutator $[\alpha^i,\beta]$ is the $i$th power of
a positive Dehn twist in the
circle $C_0$. Now let $D_j$ be the round circle in $\R^2$ with radius
$1$ centered
at the point with $x$ coordinate $3j$ and $y$ coordinate $0$. Note
that $C_0 = D_0$.
Let $\gamma$ be the composition of
a positive Dehn twist in each circle $D_j$ for $j \ge 0$, and let $\delta$
be the translation parallel to the $x$--axis through distance $3$. Then
$\gamma,\delta$
are both $C^\infty$, and the commutator $[\gamma^i,\delta]$ is the
$i$th power
of a positive Dehn
twist in the circle $D_0$. It follows that the there is an identity
$$[\alpha^i,\beta] = [\gamma^i,\delta]$$
and the group generated by $\alpha^i,\beta,\gamma^i,\delta$ is a
quotient of
$\pi_1(S)$ for $S$ the closed surface of genus $2$.

We calculate the Euler class of this action. One method is to
find a polygon $P_2$ and a developing map
$d_\rho\co \partial P_2 \to \R^2$ which is smoothed to $\delta\co \partial
P_2 \to \R^2$ with
$\text{index}(\delta) = 3+i$. It follows by
Construction~\ref{geometric_class} that
$$\rho_i^*([e])([S]) = 3 + i + 1 - 2g = i.$$
This is left as a pleasant exercise for the reader.

Alternatively, we can use Construction~\ref{algebraic_class} to
evaluate the
Euler class. Since the commutators $[\alpha^i,\beta]$ and
$[\gamma^i,\delta]$
have compact support, their images in the group of germs at infinity
are trivial.
It follows that in the group of germs, the genus $2$ surface group breaks
up into
two genus $1$ surface groups. The first $\Z \oplus \Z$ has Euler
class $i$,
as in Bestvina's example. The group generated by $\gamma^i$ and $\delta$
fixes the negative $x$--axis. It follows that we can find a section of the
universal central extension by choosing a lift which fixes a lift of the
germ at infinity of the negative $x$--axis. In particular, the extension
splits,
and the Euler class on the second $\Z \oplus \Z$ is trivial.
\end{proof}

Theorem C raises the following natural question:

\begin{qn}
Is the Euler class a bounded class on the group of real analytic
orientation-preserving diffeomorphisms of $\R^2$?
\end{qn}

Some weak evidence for a positive answer to this question is the
following:

\begin{thm}
The group $\poly^+(\R^2)$
of {\em polynomial} orientation-preserving homeomorphisms is
circularly-orderable.
\end{thm}
\begin{proof}
Let $r_0$ be the positive real axis, and let $R$ denote the set of
germs at infinity of translates of $r_0$ by orientation-preserving
polynomial homeomorphisms. Observe that $R$ contains the germs at infinity
of every straight ray. It turns out that $R$ is a naturally
circularly-ordered set.
To see this, observe that if $\alpha\co \R^2 \to \R^2$ is a polynomial
homeomorphism, then $\alpha(r_0)$ is real algebraic, and therefore its
germ at infinity is either equal to the germ at infinity of $r_0$, or is
disjoint from it. If $\lbrace r_i \rbrace$ are a finite collection
of disjoint,
properly embedded rays in $\R^2$, then they are naturally
circularly-ordered
as follows. Let $D$ be a large closed disk whose boundary
intersects the $r_i$ transversely. For each $r_i$, the intersection
$r_i \cap \partial D$ is a finite set. This set inherits a natural
ordering
from the orientation on $r_i$. For each $i$, let $p_i$ be the greatest
element
of this set. Then the rays $r_i$ are circularly-ordered by the circular
ordering
of the $p_i$ in $\partial D$. One can verify that this ordering is
independent of
the choice of (sufficiently large) $D$, and therefore defines a circular
ordering
on the germs at infinity of the $r_i$.

Now, a polynomial homeomorphism of $\R^2$ which preserves the
germ at infinity of a straight line must preserve the entire line as
a set.
Moreover, any homeomorphism of $\R^2$ which preserves {\em every}
straight line as a set must be the identity. It follows that
the natural homomorphism $\poly^+(\R^2) \to \aut(R)$
is injective. By Theorem~\ref{action_kernel}, $\poly^+(\R^2)$ is
circularly-orderable.
\end{proof}

\subsection{Torsion-free groups which are not left-orderable}

In this section we construct groups of orientation-preserving
homeomorphisms of the plane which are not left-orderable. Of course,
there
are some very simple examples of such groups: any finite cyclic group acts
on the plane by rotations, but the only finite left-orderable group
is the
trivial group. Benson Farb asked whether there are any examples
of torsion-free groups which act on the plane but are not left or
circularly-orderable. In this section, we construct
an explicit example of such a group which is not left-orderable. In the
next section, we show how to promote this example to show the existence
of a finitely presented torsion-free subgroup of $\homeo^+(\R^2)$
which is not circularly-orderable.

\begin{thm}\label{torsion_free_not_orderable}
There is a finitely presented torsion-free $C^\infty$
subgroup of $\diffeo^+(\R^2)$ which is not left-orderable.
\end{thm}
\begin{proof}
Let $G$ be given by the following presentation:
$$G = \langle a,b,c,t \; | \; a^2 = t, b^3 = t, c^7 = t, abc = t^3 \rangle.$$
Then $G$ is the fundamental group of the Seifert fibered space
$$M = (Oo0 \; | \; 3; (2,1), (3,1), (7,1) )$$
in the notation of Montesinos \cite{Montesinos}. In particular, this is
a $3$--manifold with $\til{SL(2,\R)}$ geometry, so $\til{M}$ is
homeomorphic to
$\R^3$, and $G$ is torsion-free.

On the other hand, $\pi_1(M)$ is not left-orderable, and in fact
admits no nontrivial homomorphism to $\homeo^+(\R)$. A more general fact
is proved in
\cite{Roberts_Stein_Seifert}, but it is easy enough to give a direct
proof in our case.
Suppose $G$ is left-orderable, so that we can
disjointly decompose $G = P \cup N \cup \id$ with $P\cdot P \subset P$ and
$P^{-1} = N$, by Lemma~\ref{partition}.
Without loss of generality, we can assume $t \in P$. Then since $t$
is a nontrivial
positive power of each of the generators $a,b,c$ we must have $a,b,c
\in P$ and also
$a,b,c < t$. But then $abc < t^3$ contradicting the fourth relation.

On the other hand, since $t$ is central, we can form the quotient
$\Delta$ of
$G$ by $\langle t \rangle$ with the presentation
$$\Delta = \langle a,b,c \: | \; a^2 = b^3 = c^7 = abc = \id \rangle$$
so that there is a short exact sequence
$$0 \to \Z \to G \to \Delta \to 0.$$
Now, the group $\Delta$ acts faithfully and discretely on $\R^2$ by
$C^\infty$ homeomorphisms. To see this, observe e.g. by van Kampen's
theorem
that $\Delta$ is the fundamental group of the hyperbolic $(2,3,7)$
triangle orbifold,
and therefore admits a (unique up to conjugacy and orientation)
discrete faithful representation into $\PSL(2,\R) = \isom^+(\H^2)$. So
after
choosing a diffeomorphism of $\H^2 \to \R^2$, we obtain a representation
$\rho\co  \Delta \to \diffeo^+(\R^2)$. We will construct a monotone
equivalent
(smooth) faithful action of $G$.
Since $\rho(\Delta)$ is discrete, and the set of points with nontrivial
stabilizer
in $\Delta$ are discrete, we let $O$ be the orbit of a point with
trivial stabilizer
in $\Delta$, and choose a diffeomorphism of $\R^2 \backslash O$ with
$\R^2 \backslash X$ where $X$ is the union of infinitely many
nowhere accumulating disjoint closed disks.
Let $\rho'\co G \to \diffeo^+(\R^2 \backslash X)$ be the
pushforward of the representation $\rho$ under this diffeomorphism. We
extend this
action over $X$ in the following way. Since the components of $X$ are
permuted by
$\Delta$ with trivial stabilizer, we can identify $X = \Delta \times D$
where $D$ denotes the closed unit disk. Now let the kernel $\langle t
\rangle$ of
$\rho$ act on $D$ by some smooth diffeomorphism which is fixed in a
collar neighborhood of the boundary
and has infinite order, and take the product action on $\Delta \times D$.
This product action is smooth, and since $t$ acts trivially on $\R^2
\backslash X$,
it is smooth on all of $\R^2$. This defines an extension of $\rho'$
to all of $\R^2$
which is smooth and faithful.
\end{proof}

\begin{rmk}
{\rm In \cite{Roberts_Stein_Seifert}, it is actually proved that {\em every}
action
of the fundamental group of an exceptional Seifert fibered space $M$
on $\R$ is
equivalent to the action of $\pi_1(M)$ on the space of leaves of some
taut foliation
of $M$ transverse to the circle fibration. By the Milnor--Wood inequality,
such a foliation can only exist if the Euler class of the bundle is
smaller
in absolute value than the Euler characteristic of the base orbifold,
which is
a sphere with three exceptional fibers. In particular, for any such
base orbifold, there are infinitely many Seifert fibered spaces, but
only finitely many whose fundamental groups are left-orderable. On the
other hand,
by more or less repeating the construction in
Theorem~\ref{torsion_free_not_orderable},
one can see that the fundamental group of every Seifert fibered
$3$--manifold
with infinite fundamental group admits a faithful $C^\infty$
representation in
$\diffeo^+(\R^2)$. Every such group is torsion-free. On the other hand,
all of
these groups are circularly-orderable, and contain finite index subgroups
which
are left-orderable.}
\end{rmk}

\subsection{Torsion-free groups which are not circularly-orderable}

In this subsection we show that we can combine theorem C with the
example from
Theorem~\ref{torsion_free_not_orderable} to give an example
of a torsion-free subgroup of $\homeo^+(\R^2)$ which is not
circularly-orderable,
thereby giving a complete positive answer to Farb's question. At one
point,
we make use of a powerful theorem of Tsuboi (\cite{Tsuboi_hom}) about the
homology of the group of $C^1$ homeomorphisms of $\R^n$, but aside
from this,
our construction is basically elementary.

The use of Tsuboi's theorem is not logically essential; one could avoid it
by an explicit construction. But such a construction would add
considerably
to the length of the example, without being any more enlightening.

Our proof follows the model of an argument due to Ghys \cite{Ghys_private}
which shows that
there exist planar groups which are not circularly-orderable, using
Example~\ref{Bestvina_example}, and the fact that $\homeo^+(\R^2)$
contains
torsion elements.

First we prove a couple of lemmas.

\begin{lem}\label{standard_action}
Let  $G$ be the group from Theorem~\ref{torsion_free_not_orderable}.
Then if $\rho\co G \to \homeo^+(S^1)$ is a homomorphism for which
$\rho(t)$ has
a fixed point, $\rho$ is monotone equivalent to an action which
factors through $\Delta = G/\langle t \rangle$. Moreover, $\rho$ is
either trivial, or is semi-conjugate to the ``standard'' hyperbolic
action coming
from a faithful representation of $\Delta$ in $\PSL(2,\R)$.
\end{lem}
\begin{proof}
Let $T \subset S^1$ be the fixed point set of $\rho(t)$. By blowing up
the orbit of some $p \in T$, we can assume that $T$ has nonempty interior.
Since $t$ is central, the subset $T$ is permuted by $\rho(G)$, so we
can blow
down $S^1\backslash T$ to get a monotone equivalent action on which $t$
acts trivially, so up to monotone equivalence, the action factors through
$\Delta = G/\langle t \rangle$, as claimed.

Now, it is well-known that
every homomorphism of $\Delta$ to $\homeo^+(S^1)$ is either trivial, or
monotone equivalent to a hyperbolic action.
See \cite{Calegari_forcing} for a proof, which makes use of a theorem
of Matsumoto \cite{Matsumoto}.
\end{proof}

\begin{lem}\label{diffeos_in_disk}
Let $\diffeo^+_c(\R^2)$ denote the group of diffeomorphisms of $\R^2$ with
compact support. There is a finitely-generated subgroup $L$ of
$\diffeo^+_c(\R^2)$
containing an element $t$ which is both a commutator in $L$, and
conjugate to
its square in $L$.
\end{lem}
\begin{proof}
Such subgroups are very easy to find. For example, let
$\phi\co  \R^2 \to D^2$ be a diffeomorphism such that $\|\phi'\|$
decays very
fast at infinity, and let $t$ be the pushforward by $\phi$ of a
translation. Then $t$
is certainly both a commutator, and conjugate to its square. To see that
$t$ is
a commutator, observe that the commutator of two dilations with distinct
centers
is a translation. Moreover, a dilation with dilation constant $2$
conjugates a translation to its square. If $\|\phi'\|$ decays sufficiently
quickly
(for instance, exponential decay is sufficient)
the pushforward of any element of $L$ to $D^2$ is $C^\infty$ tangent to
the identity
along $\partial D^2$.
\end{proof}

With these lemmas proved, we can now give a positive answer to Farb's
question.

\begin{thm}\label{torsion_free_not_circular}
There exists a finitely-generated, torsion-free subgroup of
$\homeo^+(\R^2)$
which is not circularly-orderable.
\end{thm}
\begin{proof}
Let $G$ be the example from Theorem~\ref{torsion_free_not_orderable}.
By abuse of notation, we will actually think of $G$ as a subgroup
of $\diffeo^+(\R^2)$. We use the same notation for the generators of $G$.
By construction, the support of $t$ is a nonaccumulating union of closed
disks $\Delta \times D$. Moreover, the action of $t$ on $1 \times D$ can
be chosen to be arbitrary, and the action on each other translate
$g \times D$ is conjugate to the action on $1 \times D$. By
Lemma~\ref{diffeos_in_disk},
we can insist that the restriction of $t$ to $1 \times D$ has the
following properties:
\begin{itemize}
\item{There are diffeomorphisms $r,s$ whose commutator satisfies
$[r,s] = t$.}
\item{There is a diffeomorphism $q$ such that $qtq^{-1} = t^2$.}
\end{itemize}
Moreover, we can choose $q,r,s$ to have support equal to
$\Delta \times D$, and further we can insist that
the action on each translate $g \times D$ is conjugate to
the action on $1 \times D$. Let $L$ denote the group generated by
$q,r,s,t$. Note that $L$ is abstractly isomorphic to the group of
the same name constructed in Lemma~\ref{diffeos_in_disk}. Define
$$G' = \langle G, L \rangle$$
Then by construction, $G'$ maps surjectively onto $\Delta$ with kernel
$L$,
and is torsion-free.

Now, if $\rho'\co G' \to \homeo^+(S^1)$ is a homomorphism, then $\rho'(t)$
is conjugate in $\homeo^+(S^1)$ to $\rho'(t^2)$. In particular, the
{\em rotation
number} of $\rho'(t)$ satisfies
$$\text{rot}(\rho'(t)) = \text{rot}(\rho'(t^2)) = 2 \text{rot}(\rho'(t))$$
and therefore this rotation number is equal to $0$, and $\rho'(t)$
has a fixed point.
For the definition and properties of rotation number, see
\cite{Calegari_forcing}.

It follows by Lemma~\ref{standard_action} that up to
monotone equivalence, $\rho'|_G$ is equivalent to the standard action of
the fundamental group of the $(2,3,7)$ hyperbolic triangle orbifold on
its circle at infinity.

The group $G$ contains a subgroup $H$ of
index $168$ corresponding to a congruence cover of the $(2,3,7)$
hyperbolic
triangle orbifold. See for instance \cite{Eightfold_way} for details.
Since $G$ is a $\Z$ extension of the orbifold group, the
subgroup $H$ is a $\Z$ extension of the fundamental
group of a closed surface $\Sigma$ of genus 3. By passing to covers in
the base
direction, for arbitrarily large $m$, we can find a
subgroup $H_m$ of $H$ of index $m$ which is a $\Z$ extension of the
fundamental
group of a closed surface $\Sigma_m$ of genus $2m+1$. The Euler class of
the extension corresponding to $H$ evaluated on the fundamental
class of the genus $3$ surface is some number $a$,
and therefore the Euler number of the extension corresponding to $H_m$
is $ma$.

So there are elements
$$\alpha_{1,m},\beta_{1,m},\dots,\alpha_{2m+1,m},\beta_{2m+1,m} \in G$$
which map
to standard generators of a surface group of genus $2m+1$ under the
homomorphism to $\Delta$, such that
$$\prod_i [\alpha_{i,m},\beta_{i,m}] = t^{ma}.$$
Now, $[q^nrq^{-n},q^nsq^{-n}]=t^{2^n}$,
so for each $m$, we can express $t^{ma}$ as a product of approximately
$\log{ma}$ commutators of this form. Together with the
$\alpha_{i,m},\beta_{i,m}$,
we can produce surface subgroups $S_m$ of $G'$ of genus approximately
$2m+2\log{ma}$. We know that the action of $\alpha_{i,m},\beta_{i,m}$ on
$S^1$ is ``standard'', up to semi-conjugacy. So by the
Milnor--Wood inequality Theorem~\ref{Milnor_Wood_inequality} and the
fact that
the restriction of $\rho'$ to $H_m$ is standard, we can estimate
the circular Euler class
$$(\rho')^*([e])[S_m] \le -4m + 4\log{ma}$$
which is as negative as we like for big $m$. Here the $-4m$ term comes
from the
standard action of $H_m$, and the $4\log{ma}$ term comes from
Milnor--Wood.
Moreover, since $r,s$ have support contained
in the union of the disks $\Delta \times D$, the formula from
Construction~\ref{geometric_class} implies that the planar Euler class of
$S_m$ is $-4m$, since we can choose a fundamental domain to calculate
the Euler class which does not intersect the support of $r,s$.

Now, by a theorem of Tsuboi,
the homology of the group of $C^1$ orientation-preserving diffeomorphisms
of $\R^2$ is equal to $\Z$ in dimension $2$; see \cite{Tsuboi} and
\cite{Tsuboi_hom}.
Let $F_2$ denote the fundamental group of a surface of genus $2$.
Let $\rho_{-4m}$ be the representation constructed in theorem C with $i
= -4m$.
Then $\rho_{-4m}(F_2)$ and $S_m$ are homologous as subgroups of the group
of $C^1$ orientation-preserving homeomorphisms of $\R^2$. By vanishing
of the oriented cobordism group in dimension $3$, we can take this
homology
to be an oriented $3$--manifold. That is, there is a $3$--manifold
$M$ with
two boundary components, $\partial(M) = \partial_1(M) \cup \partial_2(M)$,
where $\partial_1(M)$ has genus $2$ and $\partial_2(M)$ has genus
approximately
$2m+2\log ma$, and
a $C^1$ representation $\sigma_M\co  \pi_1(M) \to \homeo^+(\R^2)$ where
$\sigma_M|_{\partial_1(M)}$ is conjugate to $\rho_{-4m}|_{F_2}$, and
$\sigma_M|_{\partial_2(M)}$ is conjugate to the action of $S_m$. In fact,
$C^1$ is superfluous for our application: it suffices that each element
of $\sigma_M(\pi_1(M))$ is $C^1$ on an open dense subset. Such a
homology can
be constructed by more elementary methods, since the boundary
representation is
$C^\infty$.

We can find an oriented irreducible
$3$--manifold $N$ with torsion-free fundamental
group, and a degree $1$ map $\varphi\co N \to M$
which is a homeomorphism on the boundary. Then $\sigma_M$ induces a $C^1$
representation $\sigma_N\co \pi_1(N) \to \homeo^+(\R^2)$ with the same
boundary
behaviour. It follows that $\sigma_N$ extends to a $C^1$ action of
the amalgamated
free product
$$G'' = G'*_{S_m} \pi_1(N)$$
which is faithful on the subgroup $G'$. Let
$\sigma\co G'' \to \homeo^+(\R^2)$ denote the homomorphism defining this
action.
Note that $G''$ is torsion-free, but we do not yet know that $\sigma$ is
faithful, and therefore we don't know that $\sigma(G'')$ is torsion-free.

Let $K$ be the kernel of $\sigma$. Build a space $X$ from a $K(G',1)$ and
from $N$ by gluing along subsurfaces representing the subgroups $S_m$.
Let $\widehat{X}$ be the cover of $X$ corresponding to $K$.
Since $\sigma$ is injective on
$G'$, $\widehat{X}$ is obtained from copies of the cover $\widehat{N}$
of $N$
where $\pi_1(\widehat{N}) = K \cap \pi_1(N)$,
glued along disks in copies of the universal cover of $K(G',1)$. But
the universal
cover of a $K(\pi,1)$ is contractible, so $\widehat{X}$ is homotopic to a
$3$--manifold, obtained from copies of $\widehat{N}$ of $N$ by attaching
$1$--handles. In particular, $K$ is the fundamental group of an
irreducible
$3$--manifold. It follows that
every finitely-generated subgroup of $K$ is the fundamental group of a
compact irreducible $3$--manifold with boundary, by the Scott core
theorem~
(see \cite{Hempel}). As \cite{BRW} observed, every compact
irreducible $3$--manifold with boundary has locally indicable fundamental
group, and therefore by Theorem~\ref{locally_surjective_LO}, has LO
fundamental group. Since $K$ is locally LO, it is LO, by
Lemma~\ref{local_LO}.
We now show how to modify $\sigma$ to make it faithful.

Pick a point $p$, and blow up its orbit. That is, replace each translate
$g(p)$ of $p$ by a closed disk $D_g$ in such a way that the diameters
of these disks in
any compact region go to zero. Since $\sigma$ is $C^1$ at $p$, we can
extend the
action of $G''$ on $\R^2 \backslash G''(p)$ to the boundary circle
$\partial D_g$
by the projective linear action on the unit tangent bundle
$UT_{g(p)}\R^2$.
Let $G_p''$ be the stabilizer of $p$ in $G''$. Then $\sigma(G_p'')$
is circularly-orderable, and the stabilizer of $\partial D_g$
in $\sigma(G_p'')$
is left-orderable, by Lemma~\ref{locally_circular}. Let $K_p$ be the
stabilizer
of $\partial D_g$ in $G_p''$, so that $\sigma(K_p)$ is the stabilizer of
$\partial D_g$ in $\sigma(G_p'')$.
Since $K$ is left-orderable, and $\sigma(K_p)$ is left-orderable,
it follows that $K_p$ is left-orderable by Lemma~\ref{ses_LO}, and
$G_p''$ is circularly-orderable, by Lemma~\ref{ses_CO}. We can insert
a faithful
action of $G_p''$ on $D_g$ by thinking of $D_g$ as the cone on $S^1$, and
coning the (positive) projectivized
linear action of $G_p''$ on $S^1$ by a faithful action of $K_p$ on
$I$, which exists by Theorem~\ref{action_is_order}. We can then translate
this
action to the other disks $D_h$. The resulting action defines a
{\em faithful} representation $\sigma'$. Notice that even though
$\sigma'$ is
only $C^0$, we could not perform the blow-up without knowing that
$\sigma$
was $C^1$ at $p$. So $\sigma'$ exhibits $G''$ as a torsion-free
subgroup of
$\homeo^+(\R^2)$.

If there were a representation
$\rho''\co G'' \to \homeo^+(S^1)$ extending $\rho'$,
then the circular Euler class of $\rho_{-4m}(F_2)$ would be equal
to the circular Euler class of $S_m$, which as we calculated before,
is $\le -4m + 4\log{ma}$, violating
the Milnor--Wood inequality (Theorem~\ref{Milnor_Wood_inequality}). This
contradiction shows that $\rho''$ cannot exist, and $G''$ is not
circularly-orderable, as claimed.
\end{proof}

\begin{rmk}
{\rm Notice that in order to find a $C^\infty$ example, one just needs to find
a torsion-free $C^\infty$ homology from $S_m$ to $\rho_{-4m}(F_2)$. Such
a
homology can actually be constructed by hand,
but the construction is somewhat complicated
and unenlightening, so we do not include it here.}
\end{rmk}

\subsection{Homological rigidity for $\Z \oplus \Z$ actions}

In the remainder of the paper, we will show, in contrast to theorem C,
and Example~\ref{Bestvina_example} that for $C^1$ actions of
$\Z \oplus \Z$, the Euler class {\em always} vanishes.

We do this by studying the possible dynamics of two commuting $C^1$
homeomorphisms
$\alpha,\beta$, examining the cases based on the dynamics of $\beta$
on the
fix point set of $\alpha$.

First we treat the case that $\alpha$ has a fixed
point whose orbit under $\beta$ is not proper.

\begin{lem}\label{no_accumulation}
Let $\rho\co \Z \oplus \Z \to \homeo^+(\R^2)$ be $C^1$, and suppose that
$p$ has nontrivial stabilizer. If $p$ is fixed by a finite index subgroup,
or has
an indiscrete orbit, then $\rho^*([e])$ is trivial.
\end{lem}
\begin{proof}
Suppose the stabilizer of $p$ is a finite index subgroup. Since any finite
index subgroup of $\Z \oplus \Z$ is isomorphic to $\Z \oplus \Z$,
and since
the Euler class is multiplicative under covers, it suffices to assume that
$p$ is a global fixed point. In this case the Euler class of the planar
action is
equal to the Euler class of the (projective) action on the unit tangent
bundle at $p$ obtained by linearizing. To see this, use
Construction~\ref{geometric_class}
for a constant developing map. In particular, the Euler class is
pulled back from an action on $S^1$, and therefore by the Milnor--Wood
inequality,
the Euler class vanishes.

If $p$ is fixed by some element $\alpha$, then after passing to a finite
index subgroup, we can assume $\alpha$ is a generator. Fix some other
generator
$\beta$ and let $p_i = \beta^{i}(p)$. Suppose $\lbrace p_i \rbrace$
is indiscrete and contains
some subsequence which accumulates at $q$, which is fixed by $\alpha$.
Some subsequence approaches $q$
radially along some tangent vector $v$. Since the action is $C^1$,
it follows
that $\alpha$ fixes $q$ and $d\alpha$ fixes $v$.

Let $p_i,p_j$ be very close to $q$,
and let $\tau$ be an arc joining $p_i$ to $p_j$ so that
$\tau'(t)$ is nearly parallel to $v$ for all $t \in \tau$, and
$d\beta^{j-i}\tau'(0) = \tau'(1)$.
Then $d\alpha\tau'(0)$ is very close to $\tau'(0)$, and since $\alpha$
and $\beta$
commute, $d\alpha\tau'(1)$ is very close to $\tau'(1)$. Moreover,
near $q$,
$\alpha$ is very close in the $C^1$ topology to a linear
transformation. So, since
$d\alpha$ almost fixes $\tau'(0)$ and $\tau'(1)$, either these vectors
are almost
parallel, in which case we can take $\tau$ very close to a radial
straight arc,
or else $d\alpha$ is very close to the identity.
In either case, for our choice of $\tau$,
the arc $\alpha(\tau)$ is very close in the $C^1$ topology to $\tau$.
We pass to the finite index subgroup generated by $\alpha$ and
$\beta^{i-j}$.
Again, since the Euler class is multiplicative under covers, it
suffices to
show that it vanishes in this case.
We let $\tau \cup \alpha(\tau)$ be the image of $\partial P_1$ (where
$P_1$ is
a fundamental domain for our new copy of $\Z \oplus \Z$) under a
``degenerate'' developing map $d_\rho$ as in
Construction~\ref{geometric_class}.
Then by our geometric reasoning, the index is $1$,
and the Euler class, evaluated on the fundamental class of
$\Z \oplus \Z$, is $1 + 1 - 2 = 0$, as claimed.
\end{proof}

\subsection{The Alexander series of $\rho$}

In this section we prove that if $\alpha$ fixes some point $p$,
then if there exists $\tau$ from $p$ to $\beta(p)$ which has algebraic
intersection number zero with its translates $\beta^i(\tau)$, then the
Euler class is zero. An individual treatment of this case is not
logically necessary for the proof of theorem D, but it does make
subsequent sections easier to understand. Throughout this section and
the next, we assume that the orbit of $p$ under $\langle \beta \rangle$
is proper
(that is, it has no accumulation points) since otherwise, we could apply
Lemma~\ref{no_accumulation}.

We fix the following notation. Let $p_i = \beta^i(p)$, and let
$\tau$ be a smooth arc from $p_0$ to $p_1$ such that $d\beta(\tau'(0))
= \tau'(1)$.
Let $P = \cup_i p_i$ and $l = \cup_i \beta^i(\tau)$. We also introduce
notation
$\delta = \alpha(\tau)$.

First we introduce an algebraic tool to study the dynamics of $\alpha$ and
$\beta$, which we call the {\em Alexander series}, by analogy with
the usual
Alexander polynomial of a knot. The condition we will impose on $\tau$
will
ensure that the series is actually a Laurent polynomial.

We assume the Euler class of our representation $\rho$ is $n \ne 0$.
Look at the closed curve $\tau \cup \delta$.
This is the image of a fundamental domain for the torus under a
``degenerate'' developing map, so we can calculate the Euler class
of $\rho$
from the formula in Construction~\ref{geometric_class}.

The curves $\tau,\delta$ both represent elements in
relative homology $$[\tau],[\delta] \in H_1(\R^2,P;\Z).$$
Since there is no natural
intersection form on this group, we must work harder to find the
correct algebraic definition of the Euler class.

There is a natural Laurent series associated to the action $\rho$,
defined as follows:

\begin{defn}
{\rm The {\em Alexander series} of $\rho$, denoted
$$A_{\rho} \in \Z[[t,t^{-1}]]$$
is defined as follows.
For the moment, assume that for our choice of $\tau$, that
$\delta$ and $\tau$ are transverse at $0$.
The constant term is the {\em algebraic} intersection number, for
a generic
choice of representative (subject to the constraints above)
of $\delta$ and $\tau$, plus the algebraic sign of the intersection
at $0$.
We now define the other coefficients. For $i\ne 0$,
let the coefficient of $t^i$ be the algebraic
intersection number of the interior of $\delta$ with the interior of
$\beta^i(\tau)$.

If $d\alpha(p) = \id$, so that $\tau$ and $\delta$ cannot be made
transverse
at $0$, we let $\phi$ be a diffeomorphism centered at $p$ with small
support,
which fixes $p$ and has linear part $d\phi(p)$ a nontrivial rotation.
Then define
$$\delta' = \beta\phi\beta^{-1}\phi(\delta)$$
and compute the constant term as above for $\delta'$ and $\tau$.
We denote the $i$th coefficient of $A_{\rho}$ by $a_i$.}
\end{defn}

\begin{note}
{\rm If $\mu,\nu$ are two arcs or loops in general position, we denote their
algebraic intersection number by $\mu \cdot \nu$. So, with this notation,
we have
$$a_i = \delta \cdot \beta^i(\tau)$$
and
$$\mu \cdot \nu = - \nu \cdot \mu.$$}
\end{note}

\begin{rmk}
{\rm As defined, $A_{\rho}$ might depend on $\tau$. Nevertheless we suppress
this
in the notation.}
\end{rmk}

Notice that the coefficients of the power series only depend on the {\em
germs} of
$\tau,\delta$ at $p_0,p_1$, together with their smooth
isotopy class in $\R^2 \backslash N(p_0 \cup p_1)$, where
$N(p_0 \cup p_1)$ denotes a sufficiently small regular neighborhood of
$p_0 \cup p_1$.
So if we replace $\tau,\delta$ with
$\tau',\delta'$ which are smoothly isotopic by an isotopy supported
outside
a neighborhood of $P$, the algebraic intersection number of
$\delta'$ with $\beta^i(\tau')$ is equal to the $i$th coefficient of
$A_{\rho}$.

For the moment, we have defined $A_\rho \in \Z[[t,t^{-1}]]$. We show,
under the right conditions, that $A_\rho$ is actually a {\em
Laurent polynomial}.

But first we must see how to calculate the Euler class of $\rho$ from
the configuration
of $\tau$ and $\delta$.

\begin{lem}\label{degenerate_calculation}
Let $\rho\co \Z \oplus \Z \to \homeo^+(\R^2)$ be $C^1$. Let $\alpha,\beta$
be the
generators, and let $p$ be fixed by $\alpha$. Let $\tau$ be an embedded
arc from
$p$ to $\beta(p)$, and let $\delta$ be $\alpha(\tau)$. Let $\tau^\pm$
be immersed
proper rays in $\R^2$ such that $\tau^- \cup \tau \cup \tau^+$
is an immersed proper line, and $\tau^\pm$ intersect $\tau$ only
at their endpoints. Then the Euler class $\rho^*([e])$,
evaluated on the fundamental cycle of $\Z \oplus \Z$,
satisfies the formula
$$\rho^*([e])([\Z \oplus \Z]) = \delta \cdot \tau^+ - \delta \cdot
\tau^-.$$
\end{lem}
\begin{proof}
First we prove the lemma in the case that $\tau^+,\tau^-$ are properly
embedded,
and intersect $\tau$ only at their endpoints.

For concreteness, let $\tau^- \cup \tau \cup \tau^+$ be the real axis,
oriented
so that $\tau^-$ limits to $-\infty$.

\begin{figure}[ht]
\centerline{\relabelbox\small \epsfxsize 2.5truein
\epsfbox{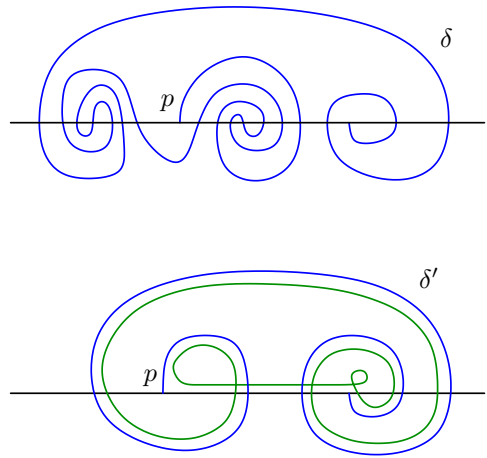}
\relabel{delta}{$\delta$}
\relabel{delta'}{$\delta'$}
\relabel{p}{$p$}
\relabel{pp}{$p$}
\endrelabelbox}
\caption{The Euler class can be calculated from the intersection numbers
$\delta \cdot \tau^\pm$. The green loop is the smoothed boundary of
a fundamental
domain for the hypothetical $\Z \oplus \Z$ action. The index of the
green loop is
$\delta \cdot \tau^+ - \delta \cdot \tau^- + 1 = 2 + 1 + 1= 4$, and
the Euler class is $3$.}
\label{degenerate_spirals}
\end{figure}

We can deform $\delta$ through a family of smooth arcs which
are all embedded, which do not intersect $p$ or $\beta(p)$ except at
their endpoints,
and whose germ near $p, \beta(p)$ is fixed. The end result of this
deformation we
will denote $\delta'$. This deformation
does not change either the index of $\delta \cup \tau$,
or the intersection numbers $\delta \cdot \tau^\pm$. After such a
deformation
we can assume $\delta'$ is of a particularly simple form: it consists
of three
pieces: a concentric spiral around $p$ through
$\delta \cdot \tau^-$ revolutions, a concentric spiral around $\beta(p)$
through $\delta \cdot \tau^+$ revolutions, and an arc between the two
which has
no vertical tangent. To see this, observe that we can deform $\delta$
to $\tau$ by
an isotopy rel. endpoints through smooth arcs,
at the cost of twisting the tangent vectors
at $p$ and at $\beta(p)$ through $\delta \cdot \tau^-$ and $\delta
\cdot \tau^+$
revolutions, respectively. Then we can ``undo'' this twisting by inserting
spirals, to produce $\delta'$.

When we smooth the curve $\tau \cup \delta$ to calculate the index by
the formula of Construction~\ref{geometric_class}, we introduce another
twist
at $\beta(p)$. As described in Construction~\ref{geometric_class}, this
is because to smooth the curve correctly, we must make a perturbation at
the degenerate edges of the developing map, in a way which is coherent
with respect to the edge pairing transformations. If the ``degenerate
edge''
from $p$ to $\alpha(p)$ is perturbed to a short vertical
edge lying on the positive side of $\tau$,
then the degenerate edge from $\beta(p)$ to $\alpha\beta(p)$ must also
be on
the positive side of $\tau$, and the end of $\delta'$ must cross over
$\tau$ near
$\beta(p)$ to reach $\alpha\beta(p)$. It follows that the index of the
smoothing of $\tau \cup \delta$ is equal to
$\delta \cdot \tau^+ - \delta \cdot \tau^- + 1$ and therefore, by the
formula of Construction~\ref{geometric_class},
the Euler class is $\delta \cdot \tau^+ -
\delta \cdot \tau^-$, as claimed.

For an example, see Figure~\ref{degenerate_spirals}.
In this case $\delta \cdot \tau^+ = 2$,
$\delta \cdot \tau^- = -1$ and the Euler class is $3$.

Note that this calculation is still valid if $\tau^-,\tau^+$ are merely
properly {\em immersed} rays which do not intersect $\tau$ except at their
endpoints. For, we can always isotope $\delta$ to $\delta'$ which is
contained
in a neighborhood of $\tau$, without affecting the intersection number
with
$\tau^+$ or $\tau^-$, and therefore the calculation is insensitive to
intersections of $\tau^\pm$ with each other or with themselves away
from $\tau$.
This completes the proof in general.
\end{proof}

With this lemma, we can establish the following fundamental properties
of $A_\rho$.

\begin{lem}\label{poly_properties}
Suppose that $\tau$ has the property that for any $i$, the
algebraic intersection number of $\tau$ and $\beta^i(\tau)$ is zero.
Then the Alexander series $A_\rho(t) = \sum_i a_i t^i$ has the following
properties:
\begin{enumerate}
\item{$A_\rho$ is actually a Laurent polynomial. That is, the coefficient
of
$t^i$ is zero for all but finitely many $i$.}
\item{The sum of the coefficients is zero:
$$A_{\rho} (1) = 0.$$}
\item{The Euler class $\rho^*([e])([\Z \oplus \Z])$, which we abbreviate
by
$e_\rho$, satisfies
$$e_\rho = \sum_{i>0} a_i - \sum_{i<0} a_i.$$}
\end{enumerate}
In this case, we call $A_\rho$ the {\em Alexander polynomial} of $\rho$.
\end{lem}
\proof
Recall our standing assumption throughout this section that
the orbit of $p$ under $\langle \beta \rangle$ is discrete, since
otherwise
Lemma~\ref{no_accumulation} would apply. We show how the first fact of
the lemma
follows from this assumption. By the assumption, for large $i$,
the points $p_i,p_{i+1}$ escape
any compact set. Since $\beta^i(\tau)$ is an
arc from $p_i$ to $p_{i+1}$, the algebraic intersection number of this arc
with $\tau \cup \delta$ is zero. But by hypothesis,
$\beta^i(\tau)$ and $\tau$ have algebraic intersection number zero,
so the algebraic intersection with $\delta$ alone is also zero.

The second fact follows for essentially the same reason as the first
fact: the
union $\tau \cup \delta$ is an (immersed) closed loop, and the sum of
coefficients is equal to the algebraic
intersection number of $\tau \cup \delta$ with the union $l=\cup_i
\beta^i(\tau)$.
Since the $p_i$ are unbounded, the algebraic intersection of
the loop and $l$ is zero, and therefore the formula follows.

The third fact follows from Lemma~\ref{degenerate_calculation}, by using
$\bigcup_{i<0} \beta^i(\tau)$ as $\tau^-$, and $\bigcup_{i>0}
\beta^i(\tau)$ as
$\tau^+$. Now, the union $\bigcup_i \beta^i(\tau)$ is {\em not}
necessarily properly immersed in
general, but note that we can homotope each $\beta^i(\tau)$ rel. a
neighborhood
of its endpoints  (not necessarily
equivariantly!) without intersecting the endpoints of $\tau$
during the track of the homotopy,
to new $\tau_i'$ until the union $\bigcup_i \tau_i'$ is properly
immersed. This does
not change the algebraic intersection number with $\delta$. Then
$$e_\rho = \sum_{i>0} \delta \cdot \tau_i' - \sum_{i<0} \delta \cdot
\tau_i'
  = \sum_{i>0} \delta \cdot \beta^i(\tau) - \sum_{i<0} \delta \cdot
\beta^i(\tau) = \sum_{i>0} a_i - \sum_{i<0} a_i.
\eqno{\smash{\raise-13pt\hbox{\qed}}}$$
%\end{proof}

With this algebraic tool, we can now prove the following lemma:

\begin{lem}\label{cancellation_lemma}
Let $\rho\co \Z \oplus \Z \to \homeo^+(\R^2)$ be a $C^1$ action which is
not free.
Suppose $\Z \oplus \Z$ is generated by $\alpha$ and $\beta$,
where $\alpha$ fixes some point $p$,
and there is a smooth embedded arc $\tau$ from $p$
to $\beta(p)$ with $d\beta(\tau'(0)) = \tau'(1)$ such that the algebraic
intersection number of $\tau$ with $\beta^i(\tau)$ is zero for any $i$.
Then the Euler class $\rho^*([e])$ is zero.
\end{lem}
\begin{proof}
Let $a_i$ be the $i$th coefficient of $A_\rho$.
Let $G_n$ be the subgroup of $\Z \oplus \Z$ of index $2n+1$ generated by
$\alpha$ and $\beta^{2n+1}$.
This gives a new representation $\rho_n$ of $\Z \oplus \Z$,
whose Euler class is equal to $2n+1$ times the Euler class of $\rho$, by
multiplicativity. The new edges $\tau_{2n+1},\delta_{2n+1}$ which are
the nondegenerate sides of a fundamental domain for the new $\Z \oplus
\Z$ subgroup
can be chosen to be the unions
$$\tau_{2n+1} = \bigcup_{i=-n}^n \beta^i(\tau),
\; \; \delta_{2n+1} = \bigcup_{i=-n}^n \beta^i(\delta)$$
respectively.

By the definition of $A_\rho$ and the fact that
$$\beta^i(\delta)\cdot \beta^j(\tau) = \delta \cdot \beta^{j-i}(\tau)
= a_{j-i}$$
we can calculate
$$A_{\rho_n}(t) = \sum_j \sum_{i=-2n}^{2n} (2n+1-|i|) a_{i+j(2n+1)} t^j$$
and so by Lemma~\ref{poly_properties}, we have
\begin{multline*}
e_{\rho_n} = (2n+1) e_\rho \\
  = \sum_{j>0} \sum_{i=-2n}^{2n} (2n+1 - |i|) a_{i+j(2n+1)}
- \sum_{j<0} \sum_{i=-2n}^{2n} (2n+1 - |i|) a_{i+j(2n+1)} \\
= \sum_{i>0} \min(i,2n+1)a_i - \sum_{i<0} \min(-i,2n+1)a_i.
\end{multline*}
Now, by Lemma~\ref{poly_properties}, the coefficients
$a_i$ are zero for $|i|$ sufficiently large. But one sees that the
coefficient of
each $a_i$ in the expression for $e_{\rho_n}$ is eventually constant
and equal to $i$,
so as $n$ increases,
the right hand side is eventually constant. On the other hand, it is
equal to
$(2n+1) e_\rho$, so $e_\rho$ is equal to zero, as claimed
\end{proof}

\subsection{The case that $\tau \cdot \beta^i(\tau) \ne 0$}

In general, we cannot find an arc $\tau$ with the properties of
Lemma~\ref{cancellation_lemma}. An example is constructed
in \cite{Handel} of an orientation-preserving homeomorphism $h$ of $\R^2$
with a proper orbit $h^i(p)$ such that every arc $\tau$ from $p$ to $h(p)$
must intersect its translates. We show how to modify our arguments to deal
with the case that $\tau \cdot \beta^i(\tau) \ne 0$ for some (possibly
infinitely many) $i$.

Let $\tau$ be any smooth arc joining $p$ to $\beta(p)$, satisfying
$d\beta(\tau'(0)) = \tau'(1)$.

\begin{defn}
{\rm The {\em self-intersection series}, denoted $B_\rho \in \Z[[t,t^{-1}]]$
is defined by
setting the constant term equal to zero, and for $i \ne 0$,
the $t^i$ term equal to the
algebraic intersection number of $\tau$ with $\beta^i(\tau)$.}
\end{defn}

There is no good way to define the homological
self-intersection of $\tau$ with itself, so
the constant term must be zero. A meaningful quantity to describe the
smooth isotopy class of $\tau$ rel. endpoints is the {\em writhe},
defined as follows:

\begin{defn}
{\rm Let $\tau$ be a smooth immersed arc in $\R^2$. Choose an
orient\-ation-preserving
identification of tangent spaces
$$\phi\co T_{\tau(0)}\R^2 \longrightarrow T_{\tau(1)}\R^2$$
such that $\phi(\tau'(0)) = \tau'(1)$.
Let $\CC$ denote the space of homotopy
classes of $C^1$ immersed curves $\sigma$ in $\R^2$ from $\tau(0)$ to
$\tau(1)$ with $\phi(\sigma'(0)) = \sigma'(1)$.
Note that these curves may intersect their endpoints in
their interiors. The class of $\tau$ in $\CC$ is denoted $[\tau]$ and
is called the
{\em writhe} of $\tau$.}
\end{defn}

Given a group $G$, an {\em affine space} for $G$ is a space $X$ with a
transitive $G$ action whose point stabilizers are trivial. For example,
the action of a group $G$ on itself makes $G$ into an affine space
for $G$.
Similarly, the fibers of a principle $G$--bundle are affine spaces
for $G$.

\begin{figure}[ht]
\centerline{\relabelbox\small \epsfxsize 2.0truein
\epsfbox{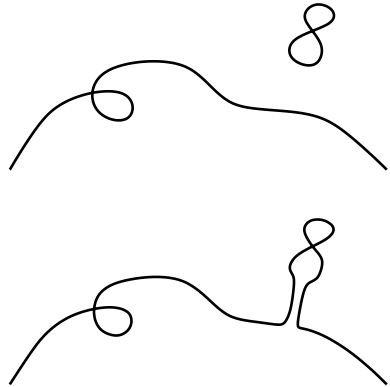}
\endrelabelbox}
\caption{The operation of connect summing with a small figure
$8$ on the positive side generates the action of $\Z$ on $\CC$.}
\label{connect_sum}
\end{figure}

Given two smooth $1$--manifolds $\kappa_1,\kappa_2$ in the plane, the
{\em connect sum} is obtained as follows. Let $I$ be an embedded
arc joining a point in the interior of $\kappa_1$ to a point in the
interior of $\kappa_2$.
Thicken $I$ to a $1$--handle $I \times I$, with boundary
$I \times \lbrace 0,1\rbrace \cup \lbrace 0,1 \rbrace \times I$, where
$I \times \lbrace 0,1 \rbrace$ is contained in the $\kappa_i$. Then
the connect sum $\kappa_1 \sharp \kappa_2$ is given by the formula
$$\kappa_1 \sharp \kappa_2 = (\kappa_1 \cup \kappa_2) \cup \partial
(I \times I)
\backslash (\kappa_1 \cup \kappa_2) \cap \partial (I \times I).$$
This produces a $1$--manifold with corners; we round the corners to get a
smooth $1$--manifold which by abuse of notation we call {\em the
connect sum}.
See Figure~\ref{connect_sum} for an illustration of this operation.

\begin{lem}
The space $\CC$ is an affine space for $\Z$ where the
action of $1$ is given by connect summing with a small
figure $8$ on the positive side.
\end{lem}
\begin{proof}
This is a very special case of Smale's immersion theorem
\cite{Smale_immersion}.
In this context, it basically follows from the Whitney trick
\cite{Whitney_trick}.
The Gauss map takes a $C^1$ immersed arc to the unit circle in $\R^2$.
As a family of arcs varies in an equivalence class of $\CC$, the
images vary
in a homotopy class rel. endpoints. The space of homotopy classes of maps
from $I$ to $S^1$ with fixed endpoints is an affine space for $\Z$; the
content of Smale's theorem in this case is that these two spaces are
equivalent.
\end{proof}

Note that there is no natural basepoint for $\CC$, but given
$\tau$ and $\delta$, the difference of the writhe of $[\tau]$ and
$[\delta]$
is an {\em integer}, which we denote by
$$[\tau] - [\delta] = n \in \Z$$
where we are using the structure of $\CC$ as an affine $\Z$ space.
That is, if $n$ is non-negative (say), then the connect sum of $\tau$
with
$n$ copies of a small figure $8$ on the positive side is in the same
equivalence
class of $\CC$ as $\delta$.

\begin{lem}\label{relative_writhe}
With $\tau,\delta$ as above, the difference in writhe is equal to the
Euler class:
$$[\tau] - [\delta] = e_\rho.$$
\end{lem}
\begin{proof}
This follows immediately from the definition of index and the
formula in Construction~\ref{geometric_class}. Firstly, observe that the
index of the smoothing only depends on the class of $\tau$ and $\delta$ in
$\CC$. If $\tau = \delta$, then the index is obviously $1$, and
$e_\rho=0$.
Finally, connect summing $\tau$ with a figure $8$ on the positive side
changes both the
writhe of $[\tau]$ and the index of the smoothing by $1$.
\end{proof}

For $\tau$ with the properties of the previous section, $B_\rho$
is identically zero.  The Alexander series $A_\rho$ is defined as
before. We now prove the analogous case of  Lemma~\ref{poly_properties}
where
$B_\rho$ is not identically zero. We prove the following
lemma for arbitrary $\tau$ joining $p$ to $\beta(p)$; in particular,
we do {\em not} assume $\tau$ is embedded.

\begin{lem}\label{poly_properties_II}
For general $\tau$, the series $A_\rho$ and $B_\rho$ satisfy the following
properties:
\begin{enumerate}
\item{The coefficients of $B_\rho$ are antisymmetric; that is, if
$B_\rho = \sum_i b_i t^i$,
then $$b_i = - b_{-i}.$$}
\item{The sum $A_\rho - B_\rho$ is a Laurent polynomial.}
\item{The sum of the coefficients of $A_\rho - B_\rho$ is zero:
$$(A_\rho - B_\rho)(1) = 0.$$}
\item{The Euler class $\rho^*([e])([\Z \oplus \Z])$, which we abbreviate
by $e_\rho$,
satisfies
$$e_\rho = \sum_{i>0} (a_i - b_i) - \sum_{i<0} (a_i - b_i).$$}
\end{enumerate}
\end{lem}
\begin{proof}
The first fact follows because the intersection number of $\tau$ with
$\beta^i(\tau)$
is equal to the intersection number of $\beta^{-i}(\tau)$ with $\tau$, and
intersection number of $1$--cycles is antisymmetric.

The second and third facts follow as in Lemma~\ref{poly_properties},
by the fact
that $\tau \cup \delta$ (where $\tau$ has the opposite orientation)
is an (immersed) closed loop, and the orbit $\beta^i(p)$ is
proper.

To calculate the Euler number, we move $\tau$ and $\delta$ by a homotopy
$\tau_t,\delta_t$
through $C^1$ immersions whose germs near $p$ and $\beta(p)$ are
constant. We move $\tau,\delta$ equivariantly, so that at each
stage, $\delta_t = \alpha(\tau_t)$. We also connect sum $\tau$ and
$\delta$ equivariantly with sufficiently many figure $8$'s on the positive
or negative sides until they are {\em embedded curves} from $p$
to $\beta(p)$.

Note that at finitely many times $t_0$, we might move $\tau_{t_0}$ through
some $p_i$, and change $b_i$ and $b_{i-1}$. However, since
$\delta_{t_0} = \alpha(\tau_{t_0})$, it follows that $\delta_{t_0}$
passes through $p_i$ simultaneously. Therefore $\tau_t \cdot
\beta^i(\tau)$
and $\delta_t \cdot \beta^i(\tau)$ change
by the same amount at such a time $t_0$, and therefore the differences
$\delta_t \cdot \beta^i(\tau) - \tau_t \cdot \beta^i(\tau)$
are always constant, and equal to $a_i - b_i$.

We also connect sum $\tau$ and $\delta$ equivariantly
with sufficiently many small figure $8$'s on positive or
negative sides. This does not affect any $a_i$ or $b_i$, and by
Lemma~\ref{relative_writhe}, it does not affect the difference in
writhe $[\tau_t] - [\delta_t]$ or the calculation of the Euler class.

It follows that the Euler class can be calculated from the
formula given in Lemma~\ref{degenerate_calculation}, after finding suitable
properly immersed rays $\tau^\pm$ such that $\tau^- \cup \tau \cup
\tau^+$ is
a properly immersed line, and $\tau^\pm$ do not intersect $\tau$ except at
their endpoints.

We obtain such rays $\tau^\pm$ by a suitable deformation of the rays
$\bigcup_{i<0} \beta^i(\tau)$ and
$\bigcup_{i>0} \beta^i(\tau)$. Initially, some $\beta^i(\tau)$ might
intersect $\tau$ essentially. But we can homotope each $\beta^i(\tau)$
through smooth curves, rel. endpoints, to $\tau_i'$ which are disjoint
from $\tau$, at the cost of crossing $\tau_i'$ $b_i$ times (algebraically)
over $p,\beta(p)$. Then for each $i$,
$$\delta \cdot \tau_i' - \tau \cdot \tau_i' = \delta \cdot \beta^i(\tau) -
\tau \cdot \beta^i(\tau).$$
Moreover, we can homotope $\tau_i'$ through immersions rel. endpoints, and
without crossing $p$ or $\beta(p)$, to some new $\tau_i''$,
so that the unions $\bigcup_{i>0} \tau_i'' = \tau^+$ and
$\bigcup_{i<0} \tau_i'' = \tau^-$ are properly immersed,
and do not intersect $\tau$ except at their endpoints. Then
$$\delta \cdot \tau_i'' = \delta \cdot \beta^i(\tau)
- \tau \cdot \beta^i(\tau) = a_i - b_i$$
so the formula for the Euler class follows.

Note that since $A_\rho - B_\rho$ is a Laurent polynomial,
this is actually a finite sum.
\end{proof}

Having established the properties of $A_\rho$ and $B_\rho$, we can now
prove the
analogue of Lemma~\ref{cancellation_lemma} for arbitrary $\tau$.

\begin{lem}\label{cancellation_lemma_II}
Let $\rho\co \Z \oplus \Z \to \homeo^+(\R^2)$ be a $C^1$ action which is
not free.
Then the Euler class $\rho^*([e])$ is zero.
\end{lem}
\begin{proof}
Without loss of generality, we can assume $\alpha$ fixes some $p$, and
the orbit of $p$ under $\beta$ is proper. Let $\tau$ be as above, and
define the series $A_\rho,B_\rho$. Let $G_n$ be the subgroup of $\Z
\oplus \Z$
of index $2n+1$ generated by $\alpha$ and $\beta^{2n+1}$ as before.
As in the proof of Lemma~\ref{cancellation_lemma}, we have
\begin{align*}
e_{\rho_n} = (2n+1)e_\rho &=
 \sum_{j>0} \sum_{i=-2n}^{2n} (2n+1-|i|)(a_{i+j(2n+1)} - b_{i+j(2n+1)}) \\
& - \sum_{j<0} \sum_{i=-2n}^{2n} (2n+1-|i|)(a_{i+j(2n+1)} - b_{i+
j(2n+1)}).
\end{align*}
On the other hand, by Lemma~\ref{poly_properties_II}, the coefficients
$a_i - b_i$ are
zero for $|i|$ sufficiently large, so as $n$ increases, the right hand
side is
eventually constant. On the other hand, this sum is equal to
$(2n+1)e_\rho$, so
$e_{\rho_n} = e_\rho = 0$, as claimed.
\end{proof}

\subsection{Vanishing of the Euler class}

From the previous section, we know that a $C^1$ $\Z \oplus \Z$ action
with nontrivial Euler class must be {\em free}; that is, no nontrivial
element has a fixed point.

We make use of several theorems, one of which was
originally proved by Brouwer (see \cite{Brouwer}),
but later given a simpler and more illuminating proof
by many people including \cite{Slaminka}, \cite{Franks}.
The first theorem says that for $\beta$ a {\em fixed-point-free}
orientation preserving homeomorphism of $\R^2$, and $p$ arbitrary, we can
find an arc $\tau$ from $p$ to $\beta(p)$ which does not intersect its
translates, except at the endpoints. More formally,

\begin{thm}[Brouwer]\label{Brouwer_theorem}
Let $\alpha\co \R^2 \to \R^2$ be orientation-preserving and fixed point
free.
Then every point has a proper orbit, and for every point $p \in \R^2$
there
is an embedded line $l \subset \R^2$ containing $p$, on which the
action of
$\alpha$ is conjugate to a translation.
\end{thm}

\begin{rmk}
{\rm The caveat of Theorem~\ref{Brouwer_theorem} is that the invariant line
$l$ it produces through every point $p$ is {\em not} generally properly
embedded, even though the orbit of every {\em point} is properly
embedded.}
\end{rmk}

We also make use of a property of fixed-point-free
orientation-preserving
homeomorphism, established by Brown, in \cite{Brown}.

Brown makes the following definition:

\begin{defn}
{\rm An orientation-preserving homeomorphism $h\co \R^2 \to \R^2$ is {\em
free} if
for every bounded, connected set $X$ with $X \cap h(X) = \emptyset$,
we have
$X \cap h^i(X) = \emptyset$ for all $i \in \Z$.}
\end{defn}

\begin{thm}[Brown]\label{freeness_theorem}
Let $h\co \R^2 \to \R^2$ be a fixed-point-free
orientation preserving homeomorphism. Then $h$ is
free.
\end{thm}

Finally, we prove an algebraic lemma:

\begin{lem}\label{equal_writhe}
Let $h\co \R^2 \to \R^2$ be a $C^1$ orientation preserving
homeomorphism. Then
for any two arcs $\tau_1,\tau_2$ from $p$ to $h(p)$
such that the union of translates of
$\tau_i$ is a $C^1$ embedded line for $i=1,2$, the writhe of $\tau_1$ and
$\tau_2$ are equal.
\end{lem}
\begin{proof}
Note that for an arc $\tau$ with the property in question, the
homeomorphism
$h$ satisfies $dh(\tau'(0)) = \tau'(1)$ so the difference in writhe of
$\tau_1$ and $\tau_2$ make sense.

We go from $\tau_1$ to $\tau_2$ by a sequence of deformations of two
different kinds.

The first kind of deformation does not change the writhe. This is a
deformation $\tau_t$ through smooth embedded curves which for each $t$
satisfy $dh(\tau'_t(0)) = \tau'_t(1)$.
For certain values of $t$, $\tau_t$ might pass through some vertex
$h^i(p)$.
By the embeddedness assumption, $i \ne 0,1$. Moreover, we choose the track
of the isotopy of $\tau_t$ to be in general position with respect to each
$h^i(p)$, so that for instance, there are only finitely many values of $t$
for which $\tau_t$ passes through some $h^i(p)$.

For any index $i$, for small values of $t-1$ we have $\tau_t \cdot
h^i(\tau_t) = 0$. As
$\tau_t$ passes over some vertex $h^i(p)$, the value of
$\tau_t \cdot h^i(\tau_t)$ changes by $e$, and the value of $\tau_t
\cdot h^{i+1}(\tau_t)$
changes by $-e$, where $e = \pm 1$ depending on
orientations. Simultaneously,
the values of $h^{-i}(\tau_t) \cdot \tau_t$ and $h^{-i-1}(\tau_t) \cdot
\tau_t$ change
by $e$ and $-e$. No other intersection numbers change for nearby values
of $t$.
It follows that the finite sum
$$w_t = \sum_{i>0} \tau_t \cdot h^i(\tau_t) - \sum_{i<0} \tau_t \cdot
h^i(\tau_t)$$
is constant under deformations of the first kind.

The second kind of deformation changes the writhe. This consists of
modifying
$\tau_t$ in a small neighborhood of $p$ to some new $\tau_{t'}$
by introducing a positive or negative ``twist'', thereby changing
the writhe
by $\pm 1$. Under such a deformation
$$\tau_{t'} \cdot h^{-1}(\tau_{t'}) - \tau_t \cdot h^{-1}(\tau_t) = 1$$
and
$$\tau_{t'} \cdot h(\tau_{t'}) - \tau_t \cdot h(\tau_t) = -1$$
for a positive twist, and the values change oppositely for a negative
twist.
It follows that $w_t$ (as defined in the previous paragraph) changes by
$\pm 2$ under
deformations of the second kind.

By combining deformations of these two kinds, we see that
$$w_t = 2 ([\tau_t] - [\tau_1])$$
where as in Lemma~\ref{relative_writhe} the notation $[\tau_t] -
[\tau_1]$ means
the difference in writhe. It follows
that $w_t$ is equal to $0$ iff the writhe of $\tau_t$ is equal to
the writhe of $\tau_1$. On the other hand, if we deform $\tau_1$ to
$\tau_2$ by
a sequence of deformations of the two kinds above, then $w_2 = 0$,
since $\tau_2$ does not intersect $h^i(\tau_2)$ except
at its endpoints. So $[\tau_2] - [\tau_1] = 0$, and the lemma is proved.
\end{proof}

\begin{rmk}
{\rm In fact, the argument of Lemma~\ref{cancellation_lemma} can be used to
give another proof of Lemma~\ref{equal_writhe}. Conversely,
Lemma~\ref{equal_writhe}
can be used to give a different proof of Lemma~\ref{cancellation_lemma}.
But the argument of Lemma~\ref{equal_writhe} does
not easily generalize to prove Lemma~\ref{cancellation_lemma_II}, whereas
the argument of Lemma~\ref{cancellation_lemma} does.}
\end{rmk}

We are now ready to prove the main theorem on homological rigidity
of $C^1$
$\Z \oplus \Z$ actions:

\begin{Z+Z_thm}\label{Z+Z_rigidity}
Let $\rho\co \Z \oplus \Z \to \homeo^+(\R^2)$ be a $C^1$ action.
Then the Euler class $\rho^*([e]) \in H^2(\Z\oplus \Z;\Z)$ is zero.
\end{Z+Z_thm}
\begin{proof}
By Lemma~\ref{cancellation_lemma_II}, we have proved this theorem except
in the
case that $\rho(\Z \oplus \Z)$ is {\em free}.

Let $p$ be arbitrary, and let $\tau$ be a smooth arc from $p$ to
$\beta(p)$ so that
$\bigcup_i \beta^i(\tau) = l$ is an embedded line. Such a $\tau$ can be
found by Theorem~\ref{Brouwer_theorem}. Let $\sigma$ be an arc from
$p$ to $\alpha(p)$. For each point $p_t \in \sigma$, we let $\tau_t$
be a smooth arc
from $p_t$ to $\beta(p_t)$ so that $\bigcup_i \beta^i(\tau_t) = l_t$ is an
embedded line. We show that for any $t_0$, there is an $\epsilon$ so that
for $|t-t_0| < \epsilon$, the family $\tau_t$ can be chosen to vary
smoothly.

To see this, observe that we can certainly choose a smoothly varying
family
$\tau_t$ near $\tau_{t_0}$ so that $\tau_t$ intersects $\beta(\tau_t)$
only
at $\beta(p_t)$. But by Theorem~\ref{freeness_theorem}, this implies that
the union of the $\beta^i(\tau_t)$ is an embedded line, as required.

It follows that we can subdivide $[0,1]$ into finitely many intervals
$$[0,q_1],[q_1,q_2],\dots,[q_n,1]$$ and choose smooth families
$\tau_t^i$ for $q_i \le t \le q_{i+1}$ with the property
that $\tau_0^0 = \tau$ and $\tau_1^n = \alpha(\tau)$. But
$\tau_{q_{i+1}}^i$ and
$\tau_{q_{i+1}}^{i+1}$ have equal writhe, by Lemma~\ref{equal_writhe},
so we can
insert a $1$--parameter family of isotopies from $\tau_{q_{i+1}}^i$ to
$\tau_{q_{i+1}}^{i+1}$ whose germ at $p_t,\beta(p_t)$ is fixed.
After inserting these isotopies, we have
constructed a one parameter family of curves $\tau_t$ from
$\tau$ to $\alpha(\tau)$ with $d\beta(\tau_t'(0)) = \tau_t'(1)$ for
all $t$. The vector field $\tau_t'$ pulls back to give a trivialization
over a fundamental domain $P_1$ for the torus.
It follows from the formula of Construction~\ref{geometric_class} that the
Euler class is zero, and the theorem is proved.
\end{proof}

\begin{rmk}
{\rm It might be possible to try to prove theorem D by a cancellation
argument, by
analogy with Lemma~\ref{cancellation_lemma} and
Lemma~\ref{cancellation_lemma_II}.
The problem is that, though the orbit of any point $p$ under any {\em
cyclic}
subgroup of a free action $\rho(\Z \oplus \Z)$ is proper, it is {\em
not} true
that the orbit under all of $\rho(\Z \oplus \Z)$ is proper.}
\end{rmk}

\begin{rmk}
{\rm Notice that the use of Brown's theorem in the proof of theorem D is
essential,
in order to be able to apply compactness of the interval $\sigma$
and compare
the relative writhe of $\tau$ and $\beta(\tau)$. Brown's theorem itself is
a kind of compactness theorem, since it lets one establish properties
of the
entire orbit of a connected set $\tau$ from the properties of only
two translates $\tau$ and $\beta(\tau)$.}
\end{rmk}

\Addresses\recd

\end{document}

%% file: gtmonout.tex
\def\ifplaintex{\expandafter\ifx\csname documentclass\endcsname\relax}

\ifplaintex 
\else
\fi

\def\gt{{\mathsurround=0pt\it $\cal G\mskip-2mu$eometry \&\ 
$\cal T\!\!$opology}}        %  journal title in recommended style

\def\gtm{{\mathsurround=0pt\it $\cal G\mskip-2mu$eometry \&\ 
$\cal T\!\!$opology $\cal M\mskip-1mu$onographs}}    %  for monographs

\def\gtp{{\mathsurround=0pt\it $\cal G\mskip-2mu$eometry \&\ 
$\cal T\!\!$opology $\cal P\!$ublications}}  % GT publications

\def\recd{{\small Received:\qua\receiveddate\ifx\reviseddate\relax
\else\qquad Revised:\qua\reviseddate\fi\par}} 

\def\volumenumber#1{\def\thevolumenumber{#1}}
\def\volumeyear#1{\def\thevolumeyear{#1}}
\def\volumename#1{\def\thevolumename{#1}}
\def\papernumber#1{\def\thepapernumber{#1}}
\def\pagenumbers#1#2{\def\startpage{#1}\def\finishpage{#2}}
\def\published#1{\def\publishdate{#1}}
\def\received#1{\def\receiveddate{#1}}
\def\revised#1{\def\reviseddate{#1}}
\def\accepted#1{\def\accepteddate{#1}}

\long\def\asciiabstract#1{\long\def\theasciiabstract{#1}}
\def\asciikeywords#1{\def\theasciikeywords{#1}}

\let\\\par
\let\thevolumenumber\relax\let\thepapernumber\relax
\let\thevolumeyear\relax\let\startpage\relax
\let\finishpage\relax\let\publishdate\relax\let\receiveddate\relax
\let\reviseddate\relax\let\accepteddate\relax\let\theasciititle\relax
\let\theasciiauthors\relax
\let\theasciiabstract\relax\let\theasciikeywords\relax

\let\theerratum\relax\let\theasciiemail\relax
\let\theshortauthors\relax\let\theshorttitle\relax

\def\startpage{1}\def\finishpage{15}\def\thepapernumber{77}

\volumenumber{2}
\volumename{Proceedings of the Kirbyfest}
\volumeyear{1999}

\long\def\maketitlep{   % start of definition of \maketitlep

\count0=\startpage

\gtm\nl        %   GT mongraphs (top left) 
{\small Volume \thevolumenumber: \thevolumename\nl 
\ifx\theerratum\relax\else Erratum \erratumnumber\nl\fi
Pages \startpage--\finishpage\nl}

\vglue 0.1truein   % top margin

{\parskip=0pt\leftskip 0pt plus 1fil\def\\{\par\smallskip}{\ifplaintex\large
\else\Large\fi\bf\thetitle}\par\medskip}   
\vglue 0.05truein 

{\parskip=0pt\leftskip 0pt plus 1fil\def\\{\par}{\sc\theauthors}
\par\medskip}%
 
\vglue 0.03truein 

{\small\leftskip 25pt\rightskip 25pt{\bf Abstract}\stdspace\theabstract

{\bf AMS Classification}\stdspace\theprimaryclass
\ifx\thesecondaryclass\relax\else; \thesecondaryclass\fi\par
{\bf Keywords}\stdspace \thekeywords\par}\vglue 7pt

}   % end of definition of \maketitlep

\font\phead=cmsl9 scaled 950
\font=cmsl9 scaled 1050
\font\pnum=cmbx10 scaled 913
\font\lnum=cmbx10 
\font\pfoot=cmsl9 scaled 950
\font=cmsl9 scaled 1050
\ifplaintex
\headline{\vbox to 0pt{\vskip -4.5mm\line{\small\phead\ifnum
\count0=\startpage ISSN 1464-8997 (on line)
1464-8989 (printed) \hfill {\pnum\folio}\else\ifodd\count0\def\\{ }% 
\ifx\theshorttitle\relax\thetitle\else\theshorttitle\fi\hfill{\pnum\folio}
\else\def\\{ and }{\pnum\folio}\hfill\ifx\theshortauthors\relax\theauthors
\else\theshortauthors\fi\fi\fi}\vss}}
\footline{\vbox to 0pt{\vglue 0mm\line{\small\pfoot\ifnum\count0=\startpage
Published \publishdate:\qua\copyright\ \gtp\hfill\else
\gtm, Volume \thevolumenumber\ (\thevolumeyear)\hfill\fi}\vss
}}
\else
\makeatletter
\def\@oddhead{{\small\lhead\ifnum\count0=\startpage ISSN 1464-8997 (on line)
1464-8989 (printed) \hfill {\lnum\number\count0}\else\ifodd\count0
\def\\{ }\ifx\theshorttitle\relax \thetitle \else\theshorttitle\fi\hfill
{\lnum\number\count0}\else\def\\{ and }{\lnum\number\count0}
\hfill\ifx\theshortauthors\relax 
\theauthors\else\theshortauthors\fi\fi\fi}}\def\@evenhead{@oddhead}
\def\@oddfoot{\small\lfoot\ifnum\count0=\startpage Published \publishdate:\qua\copyright\ \gtp\hfill\else
\gtm, Volume \thevolumenumber\ (\thevolumeyear)\hfill\fi}
\def\@evenfoot{@oddfoot}
\makeatother
\fi

\let\maketitlepage\maketitlep
\let\makeshorttitle\maketitlepage
\let\maketitle\maketitlepage

\newwrite\gtoutfile
\long\gdef\makeheadfile{  %%% start of definition of \makeheadfile
{\def\\{, }\def\s{ }
\immediate\openout\gtoutfile head.xxx
\immediate\write\gtoutfile{Proxy-for: \ifx\theasciiauthors\relax
\theauthors\else\theasciiauthors\fi\s<\ifx\theasciiemail\relax\theemail\else\theasciiemail\fi>}
\immediate\write\gtoutfile{\noexpand\\}
\immediate\write\gtoutfile{Authors: \ifx\theasciiauthors\relax
\theauthors\else\theasciiauthors\fi}
{\def\\{ }\immediate\write\gtoutfile{Title: \ifx\theasciititle\relax
\thetitle\else\theasciititle\fi}}
\immediate\write\gtoutfile{Subj-class: GT or SG, GR etc}
\immediate\write\gtoutfile{MSC-class: \theprimaryclass\ifx\thesecondaryclass\relax\else, \thesecondaryclass\fi}
\immediate\write\gtoutfile{Journal-ref: Geom. Topol. Monogr. \thevolumenumber\s
(\thevolumeyear) \startpage-\finishpage}
\immediate\write\gtoutfile{Comments: Published by Geometry and Topology Monographs at}
\immediate\write\gtoutfile{\s\s\s  http://www.maths.warwick.ac.uk/gt/GTMon\thevolumenumber/paper\thepapernumber.abs.html}
\immediate\write\gtoutfile{\noexpand\\}
\immediate\write\gtoutfile{}
\ifx\theasciiabstract\relax
\immediate\write\gtoutfile{\theabstract}\else
\immediate\write\gtoutfile{\theasciiabstract}\fi
\immediate\write\gtoutfile{}
\immediate\write\gtoutfile{\noexpand\\}
\immediate\write\gtoutfile{}
\immediate\closeout\gtoutfile}}  %%% end of definition of \makeheadfile

\def\maketitlepage{\maketitlep\makeheadfile}
\let\makeshorttitle\maketitlepage
\let\maketitle\maketitlepage